\documentclass[11pt]{article}
\newtheorem{lemma}{Lemma}[section]
\newtheorem{theorem}[lemma]{Theorem}
\newtheorem{proposition}[lemma]{Proposition}
\newtheorem{corollary}[lemma]{Corollary}
\newtheorem{definition}[lemma]{Definition}
\newtheorem{remark}[lemma]{Remark}
\newtheorem{remarks}[lemma]{Remarks}
\newtheorem{example}[lemma]{Example}

\def\square#1#2 
  {\vbox
    {\hrule
     \hbox
       {\vrule
        \vbox spread#1
          {\vfil
           \hbox spread#1{\hfil#2\hfil}
           \vfil}
        \vrule}\hrule}}
\def\sq{\square{6pt}{} }

\newenvironment{proof}{{\bf Proof}}{\hfill $ \sq $ \vskip 4mm}

\newcommand{\thlabel}[1]{\label{th:#1}}
\newcommand{\thref}[1]{Theorem~\ref{th:#1}}
\newcommand{\selabel}[1]{\label{se:#1}}
\newcommand{\seref}[1]{Section~\ref{se:#1}}
\newcommand{\lelabel}[1]{\label{le:#1}}
\newcommand{\leref}[1]{Lemma~\ref{le:#1}}
\newcommand{\prlabel}[1]{\label{pr:#1}}
\newcommand{\prref}[1]{Proposition~\ref{pr:#1}}
\newcommand{\colabel}[1]{\label{co:#1}}
\newcommand{\coref}[1]{Corollary~\ref{co:#1}}
\newcommand{\relabel}[1]{\label{re:#1}}

\newcommand{\exlabel}[1]{\label{ex:#1}}
\newcommand{\exref}[1]{Example~\ref{ex:#1}}
\newcommand{\delabel}[1]{\label{de:#1}}

\newcommand{\eqlabel}[1]{\label{eq:#1}}
\newcommand{\eqref}[1]{(\ref{eq:#1})}

\newcommand{\trr}{\triangleright}
\newcommand{\Hom}{{\rm Hom}}

\newcommand{\End}{{\rm End}}

\def\lan{\langle}
\def\ran{\rangle}

\def\im{{\rm Im\,}}

\def\trl{\triangleleft}
\newcommand{\rhu}{\mbox{$\rightharpoonup\hspace{-2ex}\rightharpoonup$}}
\newcommand{\lhu}{\mbox{$\leftharpoonup\hspace{-2ex}\leftharpoonup$}}
\def\bu{\bullet}
\def\rhua{\hbox{$\rightharpoonup$}}
\def\lhua{\hbox{$\leftharpoonup$}}

\def\pijldubbel{\lower.2ex\vbox{\hbox{$\longrightarrow$}\vspace*{-4mm}
    \hbox{$\longrightarrow$}}}
\def\doubleright#1{{\lower.2ex\vbox{\hbox{${\smash{{\cal
M}athop{\longrightarrow}
\limits^{#1}}}$}\vspace*{-4mm}\hbox{$\longrightarrow$}}}}
\def\doublerightbis#1#2{{\lower.2ex\vbox{
\hbox{${\smash{{\cal M}athop{\longrightarrow}\limits^{#1}}}$}\vspace*{-4mm}
\hbox{${\smash{{\cal M}athop{\longrightarrow}\limits_{#2}}}$}}}}
\def\doublerightleft#1#2{{\lower.2ex\vbox{
\hbox{${\smash{{\cal M}athop{\longrightarrow}\limits^{#1}}}$}\vspace*{-4mm}
\hbox{${\smash{{\cal M}athop{\longleftarrow}\limits_{#2}}}$}}}}

\def\square#1#2 
{\vbox{\hrule\hbox{\vrule\vbox spread#1
{\vfil\hbox spread#1{\hfil#2\hfil}\vfil}\vrule}\hrule}}

\def\sq{\square{5pt}{} }

\def\ot{\otimes}

\def\doublerightleft#1#2{{\lower.2ex\vbox{
\hbox{${\smash{\mathop{\longrightarrow}\limits^{#1}}}$}\vspace*{-4mm}
\hbox{${\smash{\mathop{\longleftarrow}\limits_{#2}}}$}}}}

\def\nint{\bf{Z}}

\def\klop{\hbox{$\leftharpoonup$}}

\def\ul{\underline}

\input diagrams

\begin{document}
\title{Separable functors for the category of Doi-Hopf modules. Applications}
\author{S. Caenepeel \\University of Brussels, VUB\\ Faculty of Applied
Sciences,
\\Pleinlaan 2 \\ B-1050 Brussels, Belgium
\and Bogdan Ion\thanks{Tempus visitor at UIA}\\
Department of Mathematics,\\
 Princeton University,\\
Fine Hall, Washington Road\\
Princeton, NJ 08544-1000, USA
\and G. Militaru\thanks{Research supported by the bilateral project
``Hopf algebras and co-Galois
theory" of the Flemish and Romanian governments.}
\\University of Bucharest
\\Faculty of Mathematics\\Str. Academiei 14\\RO-70109 Bucharest 1, Romania
\and Shenglin Zhu\thanks{Research supported by
the ``FWO Flanders research network WO.011.96N"}\\Institute of
Mathematics\\Fudan
University\\ Shanghai 200433, China}

\date{}
\maketitle
\begin{abstract}
\noindent We prove a Maschke type Theorem for the category of Doi-Hopf modules.
In fact, we give necessary and sufficient conditions for the functor forgetting
the $C$-coaction to be separable. This leads to a generalized notion of
integrals. Our results can be applied to obtain Maschke type Theorems for
Yetter-Drinfel'd modules, Long dimodules and modules graded by $G$-sets.
Existing Maschke type Theorems due to Doi and the authors are recovered as
special cases.
\end{abstract}

\section{Introduction}\selabel{0}
One of the key results in classical representation theory is Maschke's Theorem,
stating that a finite group algebra $kG$ is semisimple if and only if the
characteristic
of $k$ does not divide the order of $G$. Many similar results (we call them
Maschke type Theorems) exist in the literature, and the problem is always
to find
a necessary and sufficient condition for a certain object ${\cal O}$ to be
reducible or
semisimple. These objects occur in several disciplines in mathematics, they can
be for example groups, affine algebraic groups, Lie groups, locally compact
groups, or Hopf algebras, and, in a sense, this last example contains the
previous
ones as special cases. The general idea is to consider representations of
${\cal O}$.
Roughly stated, ${\cal O}$ is reducible if and only if all representations of
${\cal O}$ are completely reducible, and this comes to the fact that any
monomorphism
between two represetations splits.\\
The idea behind the proof of a Maschke type Theorem is then the following: one
applies a functor to the category of representations of ${\cal O}$, forgetting
some of the structure, and in such a way that we obtain objects in a more
handable
category, for instance vector spaces over a field or modules over a
commutative ring.
Then the strategy is to look for a deformation of a splitting map of a
monomorphism
in the ``easy" category in such a way that it becomes a splitting in the
category
of representations. In a Hopf algebraic setting, the tool that is applied
to find
such a deformation is often called an integral.\\
In fact there is more: the deformation turns out to be functorial in all cases
that are known, and this can be restated in a categorical setting: the Maschke
type Theorem comes down to the fact that the forgetful functor is {\it
separable}
in the sense of \cite{NvBvO}. A separable functor $F:\ {\cal C}\to {\cal
D}$ has more
properties: information about semisimplicity, injectivity, projectivity of
objects
in ${\cal D}$ yields information about the corresponding properties in
${\cal C}$.\\
At first sight, it seems to be a very difficult problem to decide whether a
functor is separable. However,
if a functor $F$ has a right adjoint $G$, then there
is an easy criterium for the separability of $F$: the unit $\rho:\
1\rightarrow GF$
has to be split (see \cite{Raf}), and this criterium will play a crucial
role in this
paper.\\
The classical ``Hopf algebraic Maschke Theorem" of Larson and Sweedler
(\cite{LS})
tells us that a finite dimensional Hopf algebra is semisimple if and only
if there
exists a integral that is not annihilated by the augmentation map. Several
generalizations to categories of  (generalized) Hopf modules have been
presented in the
literature, see Cohen and Fishman (\cite{CF1} and \cite{CF2}), Doi
(\cite{D0} and
\cite{D1}),
\c Stefan and Van Oystaeyen (\cite{SVO})
and the authors (\cite{CMZ3}). The results in all these papers can be
reformulated
in terms of separable functors, in fact they give sufficient conditions for a
forgetful functor to be separable. An unsatisfactory aspect is that these
conditions
are sufficient, but not necessary.\\
In this paper, we return to the setting of \cite{CMZ3}: we consider a
Doi-Hopf datum $(H,A,C)$, and the category of so-called Doi-Hopf modules
${}^C{\cal M}(H)_A$ consisting of modules with an algebra action and a coalgebra
coaction.
As explained in previous publications (\cite{CMZ1}, \cite{D2},\ldots)
${}^C{\cal M}(H)_A$
unifies modules, comodules, Sweedler's Hopf modules, Takeuchi's relative
Hopf modules,
graded modules, modules graded by $G$-sets, Long dimodules and
Yetter-Drinfel'd modules. We consider the functor $F$ forgetting the
$C$-coaction,
and we give necessary and sufficient conditions for this functor to be
separable.\\
To this end, we study natural transformations $\nu:\ GF\to 1$, and the clue
result is the following: the natural transformation $\nu$ is completely
determined if we know the map $\nu_{C\ot A}:\ C\ot C\ot A\to C\ot A$. We
recall that $C\ot A$ plays a special role in the category of Doi-Hopf modules
(although it is not a generator). Conversely, a map $C\ot C\ot A\to C\ot A$
in the category of Doi-Hopf modules can be used to construct a natural
transformation,
provided it suffices two additional properties. The next step is then to show
that the natural transformation splits the unit $\rho$ if and only if the
map $\nu_{C\ot A}$ satisfies a certain normalizing condition. We obtain a
necessary and sufficient condition for the functor $F$ to be separable, and this
condition can be restated in several ways. Actually we prove that the
$k$-algebra $V$ consisting of all natural transformations $\nu:\ GF\to 1$
is isomorphic to five different $k$-algebras, consisting of $k$-linear maps
satisfying certain properties. One of these algebras, named $V_4$,
consists of right $C^*$-linear maps $\gamma:\ C\to \Hom(C,A)$ that are
centralized by a left and right $A$-actions. We have called these maps $\gamma$
$A$-integrals, because they are closely related to Doi's total integrals
(see \seref{3.2}). The normalized $A$-integrals (also called total
$A$-integrals)
are then right units in $V_4$, and our main Theorem takes the following form:
the forgetful functor is separable if and only if there exists a total
$A$-integral
$\gamma:\ C\to \Hom(C,A)$.\\
Our technique can also be applied to the right adjoint $G$ of the functor $F$
(see \seref{2.2}), and this leads to the notion of dual $A$-integral. In
\seref{3},
we give several applications and examples: we explain how the results of
\cite{CMZ3}, \cite{D0} and \cite{D1} are special cases
(Sections \ref{se:3.2} and \ref{se:3.1}), and apply our results to some
Hopf module categories that are special cases of the Doi-Hopf module category:
Yetter-Drinfel'd modules, Long dimodules and modules graded by $G$-sets
(see Sections \ref{se:3.5}-\ref{se:3.7}). In the Yetter-Drinfel'd case, this
leads to a generalization of the Drinfel'd double in the case of an infinite
dimensional Hopf algebra, using Koppinen's generalized smash product
\cite{K} and the results of \cite{CMZ1}. Finally, if $C$ is finitely generated
and projective, then we find necessary and sufficient conditions for the
extension $A\to A\# C^*$ to be separable (\seref{3.4}). In the situation
where $C=H$,
we recover some existing results of Cohen and Fischman \cite{CF1},
Doi and Takeuchi \cite{DT}, and Van Oystaeyen, Xu and Zhang \cite{VXZ}.

\section{Preliminary results}\selabel{1}
Throughout this paper, $k$ will be a commutative ring with unit. Unless
specified otherwise, all modules, algebras, coalgebras, bialgebras,
tensor products and homomorphisms are over $k$. $H$ will be a bialgebra
over $k$, and we will extensively use Sweedler's
sigma-notation. For example, if ($C,\Delta_C$) is a coalgebra,
then for all $c\in C$ we write
$$\Delta_C(c)=\sum c_{(1)}\ot c_{(2)}\in C\ot C.$$
If ($M,\rho_M$) is a left $C$-comodule, then we write
$$\rho_M(m)=\sum m_{<-1>}\otimes m_{<0>},$$
for $m\in M$.
${}^C{{\cal M}}$ will be the category of left
$C$-comodules and $C$-colinear maps. For a $k$-algebra $A$, ${{\cal M}}_A$
(resp.
${}_A{\cal M}$) will be the category of right (resp.
left) $A$-modules and $A$-linear maps.\\
The dual $C^*=\Hom(C,k)$ of a $k$-coalgebra $C$ is a $k$-algebra. The
multiplication on $C^*$ is given by the convolution
$$\lan f*g,c\ran=\sum \lan f,c_{(1)}\ran\lan g,c_{(2)}\ran,$$
for all $f,g\in C^*$ and $c\in C$. $C$ is a $C^*$-bimodule. The
left and right action are given by the formulas
\begin{equation}\eqlabel{1.0a}
c^*\rhua c=\sum \lan c^*,c_{(2)}\ran c_{(1)}~~{\rm and}~~
c\lhua c^*= \sum \lan c^*,c_{(1)}\ran c_{(2)}
\end{equation}
for $c^*\in C^*$ and $c\in C$.
This also holds for $C$-comodules. For example if ($M,\rho_M)$ is a left
$C$-comodule, then it becomes a right $C^*$-module by
$$m\cdot c^*=\sum\lan c^*,m_{<-1>}\ran m_{<0>},$$
for all $m\in M$ and $c^*\in C^*$.
If $C$ is projective as a $k$-module, then a $k$-linear map
$f:\ M\to N$ between two left $C$-comodules is
$C$-colinear if and only if it is right $C^*$-linear.\\
An algebra $A$ that is also a left $H$-comodule is called
a left $H$-comodule algebra if the
comodule structure map $\rho_A$ is an algebra map. This means that
$$\rho_A(ab)=\sum a_{<-1>}b_{<-1>}\ot a_{<0>}b_{<0>}~~{\rm and}~~
\rho_A(1_A)=1_H\ot 1_A$$
for all $a,b\in A$.\\
Similarly, a coalgebra that is also a right $H$-module
is called a right $H$-module coalgebra if
$$\Delta_C(c\cdot h)=\sum c_{(1)}\cdot h_{(1)}\ot c_{(2)}\cdot h_{(2)}$$
and
$$\varepsilon_C(c\cdot h)=\varepsilon_C(c)\varepsilon_H(h),$$
for all $c\in C$, $h\in H$.

We recall that a functor $F:{\cal C}\to {\cal D}$ is called
{\it fully faithful} if the maps
$$\Hom_{\cal C}(M,N)\to \Hom_{\cal D}(FM,FN)$$
are
isomorphisms for all objects $M$, $N\in {\cal C}$.

\subsection{Doi-Hopf modules}\selabel{1.1}
Let $H$ be a bialgebra, $A$ a left $H$-comodule algebra and $C$ a right
$H$-module
coalgebra. We will always assume that $C$ is flat as a $k$-module.
Following \cite{CMZ1}, we will call the threetuple $(H,A,C)$
a {\it Doi-Hopf datum}. A right-left {\it Doi-Hopf module} is a $k$-module $M$
that has the structure of right $A$-module and left $C$-comodule such that the
following compatibility relation holds
\begin{equation}\eqlabel{comp}\eqlabel{1.1a}
\rho_M(m a)=\sum m_{ <-1>}\cdot a_{<-1>} \otimes m_{<0>}a_{<0>},
\end{equation}
for all $a\in A$, $m\in M$.
${}^C{\cal M}(H)_A$ will be the category of right-left Doi-Hopf modules
and $A$-linear, $C$-colinear maps.
In \cite{CR}, induction functors between
categories of Doi-Hopf modules are studied. It follows from
\cite[Theorem 1.3]{CR} that the
forgetful functor $F:\ {\cal C}={}^C{\cal M}(H)_A\to {\cal M}_A$ has a
right adjoint
$$G:\ {\cal M}_A\to {}^C{\cal M}(H)_A$$
given by
$$G(M)=C\ot M$$
with structure maps
\begin{eqnarray}
(c\ot m)\cdot a&=&\sum c\cdot a_{<-1>}\ot ma_{<0>}\eqlabel{1.1}\\
\rho_{C\ot M}(c\otimes m)&=&\sum c_{(1)}\otimes c_{(2)}\otimes m \eqlabel{1.2}
\end{eqnarray}
for any
$c\in C,~a\in A$ and $m\in M$. It is easy to see that the unit
$\rho:\ 1_{\cal C}\rightarrow GF$ of the adjoint pair $(F,G)$ is given by the
$C$-coaction $\rho_M:\ M\to M\ot C$ on any Doi-Hopf module $M$. The counit
$\delta:\ FG\rightarrow 1_{\cal C}$ is given by
$$\delta_N : C\ot N \to N,~~~\delta_N(c\ot n)=\varepsilon(c) n$$
for any right $A$-module $N$.\\
$A$ is a right $A$-module, so $G(A)=C\ot A$ and $GFG(A)=C\ot C\ot A$ are
Doi-Hopf modules. For later use, we give the action and coaction explicitely:
\begin{eqnarray}
(c\ot b)a&=& \sum ca_{<-1>}\ot ba_{<0>}\eqlabel{1.3}\\
\rho(c\ot b)&=& \sum c_{(1)}\ot c_{(2)}\ot b\eqlabel{1.4}\\
(c\ot d\ot b)a&=& \sum ca_{<-2>}\ot da_{<-1>}\ot ba_{<0>}\eqlabel{1.5}\\
\rho(c\ot d\ot b)&=& \sum c_{(1)}\ot c_{(2)}\ot d\ot b\eqlabel{1.6}
\end{eqnarray}
$\rho$ is the unit of the adjoint pair $(F,G)$, and therefore the coaction
$\rho_{C\ot A}:\ C\ot A\to C\ot C\ot A$ is $A$-linear and $C$-colinear.\\
Now assume that $H$ is a Hopf algebra. Then we can also consider the category
${}_A{\cal M}(H)^C$ of left $A$-modules that are also right $C$-comodules,
with the additional compatibility relation (see \cite{CR})
\begin{equation}\eqlabel{1.7}
\rho(am)=\sum a_{<0>}m_{<0>}\ot m_{<1>}S(a_{<-1>})
\end{equation}
Now the forgetful functor $G':\ {}_A{\cal M}(H)^C\to {\cal M}^C$ has a left
adjoint $F'$ given by $F'(M)=M\ot A$ and (see \cite{CR})
\begin{eqnarray}
a(m\ot b)&=& m\ot ab\eqlabel{1.8}\\
\rho'(m\ot b)&=& \sum m_{<0>}\ot b_{<0>}\ot m_{<1>}S(b_{<-1>})\eqlabel{1.9}
\end{eqnarray}
It follows in particular that $C\ot A=F'(C)\in {}_A{\cal M}(H)^C$, and this
makes $C\ot A$ into a left $A$-module and a right $C$-comodule.\\
The algebra $C^*$ is a left $H$-module algebra; the $H$-action
is given by the formula
$$\lan h\cdot c^*,c\ran=\lan c^*,c\cdot h\ran$$
for all $h\in H$, $c\in C$ and $c^*\in C^*$. The smash product $A\#C^*$ is
equal to $A\ot C^*$ as a $k$-module, with multiplication defined by
\begin{equation}\eqlabel{smash}
(a\#c^*)(b\#d^*)=\sum a_{<0>}b\#c^{*}* (a_{<-1>}\cdot d^*),
\end{equation}
for all $a, b\in A$, $c^*, d^*\in D^*$. Recall that we have a natural functor
$P:\ {}^C{\cal M}(H)_A\to {\cal M}_{A\# C^*}$ sending a Doi-Hopf module $M$
to itself,
with right $A\# C^*$-action given by
\begin{equation}\eqlabel{smash2}
m\cdot (a\# c^*)=\sum\lan c^*,m_{<-1>}\ran m_{<0>}a
\end{equation}
for any $m\in M$, $a\in A$ and $c^*\in C^*$. $P$ is fully faithful if $C$ is
projective as a $k$-module, and $P$ is an equivalence of categories if
$C$ is finitely generated and projective as a $k$-module (cf. \cite{D2}).\\
In \cite{K}, Koppinen introduced the following version of the
smash product: $\#(C,A)$ is equal
to $\Hom(C,A)$ as a $k$-module, with
multiplication given by the formula
\begin{equation}\eqlabel{Ko}
(f\bu g)(c)=\sum f\left(c_{(1)}\right)_{<0>}
g\left( c_{(2)}\cdot f\left( c_{(1)}\right)_{<-1>}\right)
\end{equation}
An easy computation shows that $\#(C,A)$ is an
associative algebra with unit $u_A\circ \varepsilon_C$. The natural map
$i:\ A\#C^*\to
\Hom(C,A)$ given by
$$i(a\# c^*)(c)=\lan c^*,c \ran a$$
is an algebra morphism which is an
isomorphism if $C$ is finitely generated and projective as a $k$-module.
We obtain the following diagram of algebra maps
$$
\begin{diagram}
A  & \rTo^{i_A}     & A\# C^*     & \lTo^{i_{C^*}}  & C^*  \\
   &                & \dTo^{i}    &                 &      \\
   &                &  \#(C,A)    &                 &
\end{diagram}
$$
where $i_A(a)=a\# \varepsilon$, $i_{C^*}=1_A\# c^*$ for any $a\in A$, $c^*\in
C^*$. It follows that $\Hom(C,A)$ can be viewed as
an $(A, A\# C^*)$-bimodule, where the actions are given by the restriction
of scalars.
For further use, we write down the actions explicitely
\begin{eqnarray}
(a\cdot f)(c)&=&(i(a\# \varepsilon)\bu f)(c)=
\sum a_{<0>}f\left(c\cdot a_{<-1>}\right)\eqlabel{1.10a}\\
(f\cdot a)(c)&=&f(c)a\eqlabel{1.10b}\\
(f\cdot c^*)(c)&=&\sum f\bigl(c_{(1)}\bigr)_{<0>}\lan c^*,
c_{(2)}\cdot f\bigl(c_{(1)}\bigr)_{<-1>}\ran\eqlabel{1.10c}
\end{eqnarray}
for any $f\in \Hom(C,A),~a\in A,~c^*\in C^* ~{\rm and}~ c\in C$.\\
Let $N$ be a left $A$-module. Then we can make $\Hom(C,N)$ into
a left $\#(C,A)$-module as follows
$$(\alpha\cdot g)(c)=\sum \alpha\bigl(c_{(1)}\bigr)_{<0>}
g\left( c_{(2)}\cdot \alpha\bigl( c_{(1)}\bigr)_{<-1>}\right)$$
for any $\alpha\in \#(C,A)$ and $g\in \Hom(C,N)$. By restriction of scalars,
$\Hom(C,N)$ becomes a left $A$-module.
For a right $A$-module $N$,
$\Hom(C,N)$ becomes a right $A$-module
in the usual way:
$$(f\cdot a)(c)=f(c)a$$
for any $f\in \Hom(C,N)$, $a\in A$ and $c\in C$. An easy verification shows that
$\Hom(C,N)$ is an $A$-bimodule if $N$ is an $A$-bimodule. We will apply this
to the $A$-bimodule $N= \#(C,A)$. In this case the left and right $A$-action
on $\Hom(C,N)=\Hom (C,\Hom (C,A))$ take the form
\begin{eqnarray}
(a\rhu \gamma)(c)&=&\sum a_{<0>}\cdot \gamma(c \cdot a_{<-1>})\eqlabel{1.11a}\\
(\gamma \lhu a)(c)&=&\gamma (c)\cdot a\eqlabel{1.11b}
\end{eqnarray}
for any $a\in A$, $\gamma \in \Hom(C,\#(C,A))$ and $c\in C$.\\
We have several functors
between the categories ${}^C{\cal M}(H)_A, ~{\cal M}_A, ~{\cal
M}_{\#(C,A)}$ and ${\cal M}_{A\# C^*}$.
$$
\begin{diagram}
{}^C{\cal M}(H)_A &\rTo^{~~~V~~~}       &{\cal M}_{\#(C,A)}
&\lTo^{\bu\ot_A\#(C,A)}& {\cal M}_A  \\
            &             &\dTo_R        &                      &       \\
            &             &{\cal M}_{A\#C^*}   &                      &       \\
\end{diagram}
$$
Here $R$ is the restriction of scalars functor, $V(M)=M$ with the induced
$\#(C,A)$-action
\begin{equation}\eqlabel{trl}
m\trl f=\sum m_{<0>} f(m_{<-1>})
\end{equation}
for all $m\in M$, $f\in \Hom(C,A)$. $V$ is fully faithful if $C$ is a
projective $k$-module. Also observe that $P=RV$, where $P$ is defined in
\eqref{smash2}.
We will also write $T=R(\bu\ot_A\#(C,A))$. For example,
the right $A$ and $C^*$-action on $P(C\ot A)$ are given by the formulas
$$(c\ot b)\cdot a=\sum c\cdot a_{<-1>}\ot ba_{<0>}~~{\rm and}~~
(c\ot b)\cdot c^*= c\lhua c^*\ot b$$
Moreover, with the action
$$a\cdot (c\ot b)=c\ot ab,$$
$C\ot A$ becomes an $(A,A\# C^*)$-bimodule.

\subsection{Separable functors}\selabel{1.2}
Let ${\cal C}$ and ${\cal D}$ be two categories, and $F:\ {\cal C}\to {\cal D}$
a covariant functor. Observe that we have two covariant functors
$$\Hom_{\cal C}(\bullet,\bullet):\ {\cal C}^{\rm op}\times {\cal C}\to
\ul{\rm Sets}
~~{\rm and}~~
\Hom_{\cal D}(F(\bullet),F(\bullet)):\ {\cal C}^{\rm op}\times {\cal C}\to
\ul{\rm Sets}
$$
and a natural transformation
$${\cal F}:\ \Hom_{\cal C}(\bullet,\bullet)\rightarrow
\Hom_{\cal D}(F(\bullet),F(\bullet))$$
Recall from \cite{NvBvO} that $F$ is called
{\it separable} if ${\cal F}$ splits, this means that we have a natural
transformation ${\cal P}:\ \Hom_{\cal D}(F(\bullet),F(\bullet))\rightarrow
\Hom_{\cal C}(\bullet,\bullet)$
such that ${\cal P}\circ{\cal F}$ is the identity natural transformation of
$\Hom_{\cal C}(\bullet,\bullet)$.\\
The terminology is motivated by the fact that a ring
extension $R\to S$ is separable (in the sense of \cite{DI})
or right semisimple (in the sense of \cite{HS}) if and
only if the restriction of scalars functor ${}_R{\cal M}\to {}_S{\cal M}$
is separable.\\
If the functor $F$ is separable, then we have the following version of
Maschke's Theorem (cf. \cite[Prop. 1.2]{NvBvO}): if $\alpha:\ M\to N$
in ${\cal C}$ is such that $F(\alpha)$ splits or co-splits in ${\cal D}$,
then $f$ splits or co-splits in ${\cal C}$.\\
Now suppose that $F:\ {\cal C}\to {\cal D}$ has a right adjoint $G$,
and write $\rho:\ 1_{\cal C}\rightarrow GF$ and $\delta:\ FG\rightarrow
1_{\cal D}$
for the unit and counit of this adjunction. Then we have the following
result (see \cite{Raf} and \cite{R}):

\begin{theorem}\thlabel{1.1}
Let $G:\ {\cal D}\to {\cal C}$ be a right adjoint of $F:\ {\cal C}\to {\cal
D}$.\\
1) $F$ is separable if and only if $\rho$ splits, this means that there is
a natural transformation $\nu:\ GF \rightarrow 1_{\cal C}$ such that
$\nu\circ \rho$ is the identity natural transformation of ${\cal C}$, or
$\nu_M\circ \rho_M=I_M$ for all $M\in {\cal C}$.\\
2) $G$ is separable if and only if $\delta$ co-splits, this means that there
is a natural transformation $\theta:\ FG\rightarrow 1_{\cal D}$ such that
$\delta\circ \theta$ is the identity natural transformation of ${\cal D}$.
\end{theorem}

As we already remarked, a separable functor reflects split sequences. However,
not every functor that reflects split sequences is separable: it suffices
to take a nonseparable finite field extensions $k\to l$. The restriction of
scalars functor reflects split sequences, but is not separable.

\section{Separability Theorems for Doi-Hopf modules}\selabel{2}
\subsection{The forgetful functor}\selabel{2.1}
In this Section, we will give necessary and sufficient conditions for the
forgetful functor $F:\ {}^C{\cal M}(H)_A\to {\cal M}_A$ to be
separable. This will lead to generalized notions of integral and
separability idempotent; in \seref{3}, we will discuss special cases,
and these will explain our terminology.
We will recover some existing results, and obtain some new applications.
We will always assume that $C$ is flat as a $k$-module, this is
of course no problem if we work over a field $k$.\\
In the sequel, we will write ${\cal C}={}^C{\cal M}(H)_A$, $F$ will be the
forgetful functor, and $G$ its right adjoint.
Having \thref{1.1} in mind, we will first study the natural
transformations $\nu:\ GF\rightarrow 1_{\cal C}$. Recall that such a
natural transformation consists of the following data: for every
$M\in {\cal C}$, we have an $A$-linear, $C$-colinear map
$\nu_M:\ C\ot M\to M$ satisfying the following naturality condition:
for any map $f:\ M\to N$ in ${}^C{\cal M}(H)_A$, we have a commutative
diagram
$$\begin{diagram}
C\ot M&\rTo^{\nu_M}&M\\
\dTo^{I_C\ot f}&&\dTo_{f}\\
C\ot N&\rTo^{\nu_N}&N\\
\end{diagram}$$
From \leref{2.1.2}, it will follow that the natural transformations
$\nu:\ GF\rightarrow 1_{\cal C}$ form a set $V$. Actually $V$ is a $k$-algebra;
the addition and scalar multiplication are given by the formulas
$$(\nu+\nu')_M=\nu_M+\nu'_M~~{\rm and}~~(x\nu)_M=x\nu_M$$
for all $x\in k$. The multiplication is defined as follows: for two
natural transformations $\nu$ and $\nu'$, we define
$\nu\bullet\nu'=\nu'\circ\rho\circ\nu$. In general, $V$ does not have a unit,
but if $\nu$ is a splitting of $\rho$, then $\nu$ is a right unit.
We will first give descriptions of $V$ as a $k$-module, and come back to
the algebra structure afterwards.

\begin{lemma}\lelabel{2.1.1}
Let $H$ be a Hopf algebra, and $(H,A,C)$ a Doi-Hopf datum.
With notation as above, consider a natural transformation
$\nu:\ GF\rightarrow 1_{\cal C}$. We use $\nu$ also as a notation for
$\nu=\nu_{C\ot A}:\ C\ot C\ot A\to C\ot A$.
$\nu$ then satisfies the following properties:
\begin{eqnarray}
&&\nu(c\ot d\ot ba)=b\nu(c\ot d\ot a)\eqlabel{2.1.1.1a}\\
&&\sum (d_{(2)}\ot 1_A)(\varepsilon_C\ot \rho_A)\nu(c\ot d_{(1)}\ot
1)=\nu(c\ot d\ot 1)
\eqlabel{2.1.1.1b}
\end{eqnarray}
for all $a,b\in A$ and $c,d\in C$.
\end{lemma}

\begin{proof}
For any $b\in A$, consider the map
$$f_b:\ C\ot A\to C\ot A,~~f_b(c\ot a)=b(c\ot a)=c\ot ba$$
for all $a,b\in A$ and $c\in C$ (see \eqref{1.8}). It is easy to check
that $f_b$ is a morphism in ${\cal C}$, and,
from the naturality of $\nu$, it follows that we have a commutative
diagram
$$\begin{diagram}
C\ot C\ot A&\rTo^{\nu}&C\ot A\\
\dTo^{I_C\ot f_b}&&\dTo_{f_b}\\
C\ot C\ot A&\rTo^{\nu}&C\ot A
\end{diagram}$$
and \eqref{2.1.1.1a} is equivalent to the commutativity of this diagram.\\
If $M$ is a $k$-module, then $C\ot A\ot M$ can be made into a Doi-Hopf
module in a natural way, by using the $A$-action and $C$-coaction on
$C\ot A$. From the naturality of $\nu$, it follows that
\begin{equation}\eqlabel{2.1.1.1extra}
\nu_{C\ot A\ot M}=\nu\ot I_M
\end{equation}
Indeed, for every $m\in M$, the map $g_m:\ C\ot A\to C\ot A\ot M$ given
by
$$g_m(c\ot a)=c\ot a\ot m$$
is in the category of Doi-Hopf modules ${\cal C}$, and we have a
commutative diagram
$$\begin{diagram}
C\ot C\ot A&\rTo^{\nu}&C\ot A\\
\dTo^{I_C\ot g_m}&&\dTo_{g_m}\\
C\ot C\ot A\ot M&\rTo^{\nu_{C\ot A\ot M}}&C\ot A\ot M
\end{diagram}$$
and we find \eqref{2.1.1.1extra}.\\
Now let $\ul{C}=C$ considered only as a $k$-module. Then the right
$C$-coaction $\rho':\ C\ot A\to C\ot A\ot C$ (see \eqref{1.9}) is a morphism in
${\cal C}$, and, using the naturality of $\nu$ and \eqref{2.1.1.1extra}, we find
a commutative diagram
$$\begin{diagram}
C\ot C\ot A&\rTo^{\nu}&C\ot A\\
\dTo^{I_C\ot \rho'}&&\dTo_{\rho'}\\
C\ot C\ot A\ot \ul{C}&\rTo^{\nu\ot I_C}&C\ot A\ot \ul{C}
\end{diagram}$$
Write $\nu(c\ot d\ot 1_A)=\sum_i c_i\ot a_i$, and apply the diagram
to $c\ot d\ot 1_A$, with $c,d\in C$. We obtain
$$\sum \nu(c\ot d_{(1)}\ot 1_A)\ot d_{(2)}=
\sum c_{i_{(1)}}\ot a_{i_{<0>}}\ot c_{i_{(2)}}S(a_{i_{<-1>}})$$
Applying $\varepsilon_C\ot I_A\ot I_C$ to both sides, and then switching
the two factors, we find
$$\sum d_{(2)}\ot (\varepsilon_C\ot I_A)\bigl(\nu(c\ot d_{(1)}\ot 1_A)\bigr)=
\sum c_{i}S\bigl(a_{i_{<-1>}}\bigr)\ot a_{i_{<0>}}$$
Now we apply $\rho_A$ to the second factor of both sides. Using the
fact that $\rho_A\circ (\varepsilon_C\ot I_A)=\varepsilon_C\ot \rho_A$,
we obtain
$$\sum d_{(2)}\ot (\varepsilon_C\ot \rho_A)\bigl(\nu(c\ot d_{(1)}\ot 1_A)\bigr)=
\sum c_{i}S\bigl(a_{i_{<-2>}}\bigr)\ot a_{i_{<-1>}}\ot a_{i_{<0>}}$$
We now find \eqref{2.1.1.1b} after we let the second factor act on the
first one.
\end{proof}

Observe that \eqref{2.1.1.1a} means that $\nu$ is left $A$-linear.
>From \eqref{2.1.1.1a}, \eqref{2.1.1.1b} and \eqref{1.3}, it follows that
\begin{equation}
\nu(c\ot d\ot a)= a\nu(c\ot d\ot 1)
= \sum (d_{(2)}\ot a)\bigl((\varepsilon_C\ot I_A)\nu(c\ot d_{(1)}\ot 1)\bigr)
\eqlabel{2.1.1.1c}
\end{equation}

We will now prove that the natural transformation $\nu$ is completely
determined by $\nu_{C\ot A}$.

\begin{lemma}\lelabel{2.1.2}
With notations as in \leref{2.1.1}, we have
\begin{equation}\eqlabel{2.1.2.1}
\nu_M(c\ot m)=
\sum m_{<0>}\bigl((\varepsilon_C\ot I_A)\nu(c\ot m_{<-1>}\ot 1)\bigr)
\end{equation}
for every Doi-Hopf module $M$, $m\in M$ and $c\in C$.
\end{lemma}

\begin{proof}
First take $N=C\ot M$, where $M$ is a right $A$-module. For any $m\in M$,
the map $g_m:\ C\ot A\to C\ot M$, $g(c\ot a)=c\ot ma$ is a morphism in
${\cal C}$, so we have a commutative diagram
$$\begin{diagram}
C\ot C\ot A&\rTo^{\nu_{C\ot A}}&C\ot A\\
\dTo^{I_C\ot g_m}&&\dTo_{g_m}\\
C\ot C\ot A&\rTo^{\nu_{C\ot M}}&C\ot A
\end{diagram}$$
Evaluating this diagram at $c\ot d\ot 1_A$, we obtain
\begin{eqnarray}
\nu_{C\ot M}(c\ot d\ot m)&=& g_m(\nu_{C\ot A}(c\ot d\ot 1_A))\nonumber\\
&=& \sum (d_{(2)}\ot m)\bigl((\varepsilon_C\ot I_A)\nu(c\ot d_{(1)}\ot
1_A)\bigr)
\eqlabel{2.1.2.1a}
\end{eqnarray}
Now consider an arbitrary Doi-Hopf module $M\in {\cal C}$. The coaction
$\rho_M:\ C\to C\ot M$ is a morphism in ${\cal C}$, and we have a
commutative diagram
$$\begin{diagram}
C\ot M&\rTo^{\nu_{M}}&M\\
\dTo^{I_C\ot \rho_M}&&\dTo_{\rho_M}\\
C\ot C\ot M&\rTo^{\nu_{C\ot M}}&C\ot M
\end{diagram}$$
We apply the diagram to $c\ot m$. Using \eqref{2.1.2.1a}, we obtain
\begin{eqnarray*}
\rho_M(\nu_M(c\ot m))&=&
\nu_{C\ot M}(\sum c\ot m_{<-1>}\ot m_{<0>})\\
&=& \sum (m_{<-1>}\ot m_{<0>})((\varepsilon_C\ot I_A)\nu(c\ot m_{<-2>}\ot 1_A))
\end{eqnarray*}
and \eqref{2.1.2.1} follows after we apply $\varepsilon_C\ot I_M$ to
both sides.
\end{proof}

Now let $V_1$ be the $k$-module consisting of all maps
$\nu\in\Hom_{\cal C}(C\ot C\ot A,C\ot A)$ satisfying \eqref{2.1.1.1a}
and \eqref{2.1.1.1b}. $\nu\in V_1$ is called {\it normalized} or {\it
splitting} if
\begin{equation}\eqlabel{2.1.1.1d}
\sum \nu(c_{(1)}\ot c_{(2)}\ot a)=c\ot a
\end{equation}
for all $c\in C$ and $a\in A$. Recall that $V$ is the $k$-module consisting of
all natural transformations $\nu:\ GF\to 1_{\cal C}$. With these notations, we
have the following result:

\begin{theorem}\thlabel{2.1.3}
Let $H$ be a  Hopf algebra, and $(H,A,C)$ a Doi-Hopf datum.
Then the $k$-modules $V$ and $V_1$ are isomorphic, and
normalized elements in $V_1$ correspond to splittings of the unit
$\rho:\ 1_{\cal C}\to GF$.\\
Consequently $F$ is a separable functor if and only if there exists a
normalized element in $V_1$.
\end{theorem}

\begin{proof}
For $\nu\in V_1$, we define $g(\nu)=\nu$ as follows: for every $M\in {\cal C}$,
$\nu_M:\ C\ot M\to M$ is given by
\begin{equation}\eqlabel{2.1.3.1a}
\nu_M(c\ot m)=\sum m_{<0>}((\varepsilon_C\ot I_A)\nu(c\ot m_{<-1>}\ot 1_A))
\end{equation}
(cf. \leref{2.1.2}). We have to show that $\nu$ is indeed a natural
transformation. It suffices to show that $\nu_M$ is right $A$-linear and
left $C$-colinear, and that the naturality condition is satisfied.\\
a) $\nu_M$ is right $A$-linear. For all $m\in M,~c\in C$ and $a\in A$, we have
\begin{eqnarray*}
\nu_M(c\ot m) a&=&
\sum m_{<0>}\bigl((\varepsilon_C\ot I_A)(\nu(c\ot m_{<-1>}\ot 1))\bigr)a\\
&=& \sum m_{<0>}\Bigl((\varepsilon_C\ot I_A)\nu((c\ot m_{<-1>}\ot 1))a\bigr)\\
&=& \sum m_{<0>}\Bigl((\varepsilon_C\ot I_A)\bigl(
\nu(c\cdot a_{<-2>}\ot m_{<-1>}\cdot a_{<-1>}\ot a_{<0>})\bigr)\Bigr)\\
\eqref{2.1.1.1a}~~~&=& \sum (m_{<0>}a_{<0>})\Bigr((\varepsilon_C\ot I_A)
\bigl(
\nu(c\cdot a_{<-2>}\ot m_{<-1>}\cdot a_{<-1>}\ot 1)\bigr)\Bigr)\\
&=& \nu_M(\sum c\cdot a_{<-1>}\ot ma_{<0>})
\end{eqnarray*}
where we used the fact that
$\varepsilon_C\ot I_A:\ C\ot A\to A$ is left and right $A$-linear.\\
b) $\nu_M$ is left $C$-colinear. Write
$$\nu(c\ot d\ot 1)=\sum_i c_i\ot a_i$$
>From the fact that $\nu$ is left $C$-colinear, it follows that
$$\sum c_{(1)}\ot \nu(c_{(2)}\ot d\ot 1)=\rho(\nu(c\ot d\ot 1))=
\sum_i c_{i_{(1)}}\ot c_{i_{(2)}}\ot a_i$$
Applying $I_C\ot \varepsilon_C\ot I_A$ to both sides, we obtain
\begin{equation}\eqlabel{2.1.3.1}
\sum c_{(1)}\ot (\varepsilon_C\ot I_A)(\nu(c_{(2)}\ot d\ot 1))=
\sum_i c_i\ot a_i=\nu(c\ot d\ot 1)
\end{equation}
We temporarily introduce the following notation, for $m\in M,~c\in C$ and
$a\in A$:
$$(I\ot m)(c\ot a)=c\ot ma$$
Now
\begin{eqnarray*}
\rho_M(\nu_M(c\ot m))&=&
\rho_M\Bigl(\sum m_{<0>}\bigl((\varepsilon_C\ot I_A)
\nu(c\ot m_{<-1>}\ot 1_A)\bigr)\Bigr)\\
&=& \sum \rho_M(m_{<0>})\rho_A\Bigl((\varepsilon_C\ot I_A)
\nu(c\ot m_{<-1>}\ot 1_A)\Bigr)\\
&=& \sum (m_{<-1>}\ot m_{<0>})\rho_A\Bigl((\varepsilon_C\ot I_A)
\nu(c\ot m_{<-2>}\ot 1_A)\Bigr)\\
\eqref{2.1.1.1b}~~~&=&
\sum (I\ot m_{<0>})\nu(c\ot m_{<-1>}\ot 1_A)\\
\eqref{2.1.3.1}~~~&=&
\sum (c_{(1)}\ot m_{<0>})\cdot \Bigl((\varepsilon_C\ot I_A)
\bigl(\nu(c_{(2)}\ot m_{<-1>}\ot 1_A)\bigr)\Bigr)\\
&=& \sum c_{(1)}\ot \nu_M(c_{(2)}\ot m)\\
&=& (I_C\ot \nu_M)(\rho_M(c\ot m))
\end{eqnarray*}
c) We now prove the naturality condition: let $\alpha:\ M\to N$ is
a morphism in ${\cal C}$. Then
\begin{eqnarray*}
\nu_N(c\ot\alpha(m))&=&
\sum \alpha(m)_{<0>}\bigl((\varepsilon_C\ot I_A)
\nu(c\ot \alpha(m)_{<-1>}\ot 1)\bigr)\\
(\alpha~{\rm is~left~}C{\rm~colinear})~~
&=& \sum \alpha(m_{<0>})\bigl((\varepsilon_C\ot I_A)
\nu(c\ot m_{<-1>}\ot 1)\bigr)\\
(\alpha~{\rm is~right~}A{\rm~linear})~~
&=& \sum \alpha\Bigl(m_{<0>}\bigl((\varepsilon_C\ot I_A)
\nu(c\ot m_{<-1>}\ot 1)\bigr)\Bigr)\\
&=& \alpha(\nu_M(c\ot m))
\end{eqnarray*}
Now define $f:\ V\to V_1$ by
$$f(\nu)=\nu_{C\ot A}$$
>From \leref{2.1.1}, it follows that $f(\nu)\in V_1$, and from
\leref{2.1.2} that $g\circ f=I_V$.
>From \eqref{2.1.1.1c}, it follows that $f\circ g=I_{V_1}$.\\
Finally, if $\nu$ is normalized, then for $M\in {\cal C}$ and $m\in M$,
we have
\begin{eqnarray*}
\nu_M(\sum m_{<-1>}\ot m_{<0>})&=&
\sum m_{<0>}\bigl((\varepsilon_C\ot I_A)\nu(m_{<-2>}\ot m_{<-1>}\ot 1_A)\bigr)\\
&=& \sum m_{<0>}\bigl((\varepsilon_C\ot I_A)(m_{<-1>}\ot 1_A)\bigr)=m
\end{eqnarray*}
and $\rho$ is split by $\nu$. The other statements of the Theorem are
immediate consequences.
\end{proof}

\begin{remark}\relabel{2.1.3a}\rm
Observe that the antipode $S$ of $H$ is used in the construction of
$f:\ V\to V_1$, but not in the construction of $g:\ V_1\to V$. So if
$H$ is only a bialgebra, then we are still able to define a homomorphism
$g:\ V_1\to V$, and we can still conclude that the existence of a
normalized element in $V_1$ is a sufficient condition for the separability
of the forgetful functor $F$. We will come back to this situation in
\seref{2.3}.
\end{remark}
The $k$-module $V_1$ can be rewritten in at least four different ways.
Consider the following $k$-modules $V_2,~V_3,~V_4$ and $V_5$.

{\bf 2)} $V_2$ consists of all maps $\lambda:\ C\ot C\to C\ot A$
satisfying the following properties:
\begin{eqnarray}
&&\sum (d_{(2)}\ot 1_A) (\varepsilon \ot \rho_A)\lambda (c\ot d_{(1)})=
\lambda (c\ot d)\eqlabel{2.1.1.2a}\\
&&\sum a_{<0>} \lambda (c\cdot a_{<-2>}\ot d\cdot a_{<-1>})=
\lambda (c\ot d)a
\eqlabel{2.1.1.2b}\\
&&\sum c_{(1)} \ot \lambda (c_{(2)}\ot d)=\rho_{C\ot A}(\lambda (c\ot d))
\eqlabel{2.1.1.2c}
\end{eqnarray}
for all $c,d\in C$. $\lambda\in V_2$ is called normalized if
\begin{equation}\eqlabel{2.1.1.2d}
\sum \lambda(c_{(1)}\ot c_{(2)})=c\ot 1
\end{equation}
for all $c\in C$.

{\bf 3)} $V_3$ consists of all maps $\theta:\ C\ot C\to A$
satisfying the following properties:
\begin{eqnarray}
&&\theta (c\ot d) a= \sum a_{<0>} \theta (c\cdot a_{<-2>}\ot d\cdot a_{<-1>})
\eqlabel{2.1.1.3a}\\
&&\sum c_{(1)}\ot \theta(c_{(2)}\ot d)=\sum d_{(2)}\theta(c\ot d_{(1)})_{<-1>}
\ot \theta(c\ot d_{(1)})_{<0>}\eqlabel{2.1.1.3b}
\end{eqnarray}
for all $c,d\in C$ and $a\in A$. $\theta\in V_3$ is called normalized,
or a {\it coseparability idempotent} if
\begin{equation}\eqlabel{2.1.1.3c}
\sum\theta(c_{(1)}\ot c_{(2)})=\eta_A(\varepsilon_C(c))
\end{equation}
for all $c\in C$ (cf. \exref{2.1.9}).

{\bf 4)} $V_4$ consists of all maps $\gamma:\ C\to\Hom(C,A)$
satisfying the following properties:\\
a) $\gamma$ is centralized by the left and right actions of $A$ on
$\Hom(C,A)$:
\begin{equation}\eqlabel{2.1.1.4a}
a\rhu \gamma=\gamma \lhu a
\end{equation}
for any $a\in A$, or, equivalently,
\begin{equation}\eqlabel{2.1.1.4b}
\sum a_{<0>}\cdot\gamma(c\cdot a_{<-1>})=\gamma(c)\cdot a
\end{equation}
for all $a\in A$ and $c\in C$, or
\begin{equation}\eqlabel{2.1.1.4c}
\sum a_{<0>}\gamma(c\cdot a_{<-2>})(da_{<-1>})=(\gamma(c)(d)) a
\end{equation}
for all $a\in A$ and $c,d\in C$.\\
b) For all $c,d\in C$, we have
\begin{equation}\eqlabel{2.1.1.4e}
\sum c_{(1)}\ot \gamma(c_{(2)})(d)=
\sum d_{(2)}\gamma(c)(d_{(1)})_{<-1>}\ot \gamma(c)(d_{(1)})_{<0>}
\end{equation}
\eqref{2.1.1.4e} looks artificial, but it has a natural interpretation in
the case where $C$ is projective as a $k$-module. Then \eqref{2.1.1.4e}
is equivalent to
\begin{equation}\eqlabel{2.1.1.4d}
\sum\lan c^*,c_{(1)}\ran\gamma(c_{(2)})(d)=
\sum\lan c^*,d_{(2)}\gamma(c)(d_{(1)})_{<-1>}\ran\gamma(c)(d_{(1)})_{<0>}
\end{equation}
which means that $\gamma$ is right $C^*$-linear, where the right
$C^*$-actions on
$C$ and $\Hom(C,A)$ are given by the formulas \eqref{1.0a} and \eqref{1.10c}.\\
An element of $V_4$ is called normalized if
\begin{equation}\eqlabel{2.1.1.4f}
\sum \gamma(c_{(1)})(c_{(2)})=\varepsilon(c)1_A
\end{equation}
for all $c\in C$.\\
{\bf 5)} $V_5$ consists of all $A$-bimodule maps $\psi:\ C\ot A\to\Hom(C,A)$
satisfying the additional condition
\begin{equation}\eqlabel{2.1.1.5a}
\sum c_{(1)}\ot \psi(c_{(2)}\ot 1_A)(d)=
\sum d_{(2)}\psi(c\ot 1_A)(d_{(1)})_{<-1>}\ot \gamma(c)(d_{(1)})_{<0>}
\end{equation}
for all $c,d\in C$. If $C$ is projective as a $k$-module, then a map
$\psi\in V_5$
is nothing else then an $(A,A\# C^*)$-bimodule map.\\
$\psi\in V_5$ is called normalized if
\begin{equation}\eqlabel{2.1.1.5}
\sum \psi(c_{(1)}\ot 1_A)(c_{(2)})=\varepsilon(c)1_A
\end{equation}
for all $c\in C$.

\begin{proposition}\prlabel{2.1.4}
Let $H$ be a bialgebra, and $(H,A,C)$ a Doi-Hopf datum. Then the $k$-modules
$V_1,\ldots,V_5$ are isomorphic, and normalized elements correspond
to normalized elements.
\end{proposition}

\begin{proof}
1) We define $f_1:\ V_1\to V_2$ and $g_1:\ V_2\to V_1$ by
$$f_1(\nu)=\lambda~~{\rm and}~~g_1(\lambda)=\nu$$
with
\begin{equation}\eqlabel{2.1.4.1}
\lambda(c\ot d)=\nu(c\ot d\ot 1)~~;~~
\nu(c\ot d\ot a)=a\lambda(c\ot d)
\end{equation}
It is straightforward to see that \eqref{2.1.1.2b} is equivalent to
right $A$-linearity of $\nu$, \eqref{2.1.1.2c} to left $C$-colinearity of
$\nu$, and \eqref{2.1.1.1b} to \eqref{2.1.1.2a}. \eqref{2.1.1.1a}
follows from \eqref{2.1.4.1}. It is also clear that $f_1$ and $g_1$ are
each others inverses.\\
2) $f_2:\ V_2\to V_3$ and $g_2:\ V_3\to V_2$ are defined by
$$f_2(\lambda)=\theta~~{\rm and}~~g_2(\theta)=\lambda$$
with
\begin{eqnarray*}
\theta(c\ot d)&=& (\varepsilon_C\ot I_A)(\lambda(c\ot d))\\
\lambda(c\ot d)&=& \sum c_{(1)}\ot\theta(c_{(2)}\ot d)
\end{eqnarray*}
We first show that $f_2$ is well-defined: if $\lambda\in V_2$, then
$\theta=f_2(\lambda)$ satisfies \eqref{2.1.1.3a} and \eqref{2.1.1.3b}.
We apply $\varepsilon_C\ot I_A$ to both sides of \eqref{2.1.1.2b}.
Writing $\lambda(c\ot d)=\sum_i c_i\ot a_i$, we obtain for the
right hand side:
\begin{eqnarray*}
(\varepsilon_C\ot I_A)(\lambda(c\ot d)a)&=&
(\varepsilon_C\ot I_A)(\sum_i c_ia_{<-1>}\ot a_ia_{<0>})\\
&=& \sum_i \varepsilon_C(c_i)a_ia=\theta(c\ot d)a
\end{eqnarray*}
For all $a\in A$, $c,d\in C$, with again $\lambda(c\ot d)=\sum_i c_i\ot a_i$,
we also have
$$(\varepsilon_C\ot I_A)(a\lambda(c\ot d))=\sum \varepsilon(c_i)aa_i=
a\theta(c\ot d)$$
so the left hand side amounts to
$$\sum a_{<0>}\theta(ca_{<-2>}\ot da_{<-1>})$$
and \eqref{2.1.1.3a} follows.\\
Now apply $I_C\ot \varepsilon_C\ot I_A$ to both sides of \eqref{2.1.1.2c}.
This gives
\begin{eqnarray*}
\sum c_{(1)}\ot \theta(c_{(2)}\ot d)&=&
((I_C\ot \varepsilon_C\ot I_A)\circ\rho_{C\ot A}\circ \lambda)(c\ot d)\\
&=& \lambda(c\ot d)\\
\eqref{2.1.1.2a}~~~
&=& \sum d_{(2)}(\varepsilon_C\ot \rho_A)\lambda(c\ot d_{(1)})\\
&=& \sum d_{(2)}(\theta(c\ot d_{(1)}))_{<-1>}\ot (\theta(c\ot d_{(1)}))_{<0>}
\end{eqnarray*}
proving \eqref{2.1.1.3b}.\\
Conversely, if $\theta\in V_3$, then $\lambda=g_2(\theta)\in V_2$. Indeed,
\begin{eqnarray*}
\sum d_{(2)}(\varepsilon_C\ot \rho_A)\lambda(c\ot d_{(1)})&=&
\sum d_{(2)}\rho_A(\theta(c\ot d_{(1)}))\\
&=& \sum d_{(2)}(\theta(c\ot d_{(1)}))_{<-1>}\ot (\theta(c\ot d_{(1)}))_{<0>}\\
\eqref{2.1.1.3b}~~~&=&
\sum c_{(1)}\ot \theta(c_{(2)}\ot d)=\lambda(c\ot d)
\end{eqnarray*}
proving \eqref{2.1.1.2a}. Furthermore
\begin{eqnarray*}
\sum a_{<0>}\lambda(ca_{<-2>}\ot da_{<-1>})&=&
\sum c_{(1)}a_{<-3>}\ot a_{<0>}\theta(c_{(2)}a_{<-2>}\ot da_{<-1>})\\
\eqref{2.1.1.3a}~~~&=&
\sum c_{(1)}a_{<-1>}\ot\theta(c_{(2)}\ot d)a_{<0>}\\
&=& \lambda(c\ot d)a
\end{eqnarray*}
proving \eqref{2.1.1.2b}. Finally
\begin{eqnarray*}
\sum c_{(1)}\ot \lambda(c_{(2)}\ot d)&=&
\sum c_{(1)}\ot c_{(2)}\ot \theta(c_{(3)}\ot d)\\
&=& \rho_{C\ot A}(\lambda(c\ot d))
\end{eqnarray*}
proving \eqref{2.1.1.2c}. It is clear that $g_3$ and $f_3$ are each others
inverses.\\
3) $f_3:\ V_3\to V_4$ and $g_3:\ V_4\to V_3$ are defined by
$$f_3(\theta)=\gamma~~{\rm and}~~g_3(\gamma)=\theta$$
with
$$\gamma(c)(d)=\theta(c\ot d)$$
It can be checked easily that \eqref{2.1.1.3a} is equivalent to
\eqref{2.1.1.4a}, and \eqref{2.1.1.3b} to \eqref{2.1.1.4e}, so
$f_3$ and $g_3$ are well-defined. Obviously $g_3=f_3^{-1}$.\\
4) $f_4:\ V_4\to V_5$ and $g_4:\ V_5\to V_4$ are defined by
$$f_4(\gamma)=\psi~~{\rm and}~~g_4(\psi)=\gamma$$
with
$$\psi(c\ot a)=a\cdot \gamma(c)~~{\rm and}~~
\gamma(c)=\psi(c\ot 1_A)$$
Let us show that $f_4$ is well-defined: $f_4(\gamma)=\psi\in V_5$ if
$\gamma\in V_4$. $\psi$ is left $A$-linear since
\begin{eqnarray*}
\psi(a\cdot(c\ot b))&=&\psi(c\ot ab)\\
&=&ab\cdot \gamma(c)\\
&=&a\cdot(\psi(c\ot b)),
\end{eqnarray*}
right $A$-linear since
\begin{eqnarray*}
\psi((c\ot b)\cdot a)&=&\psi(\sum c\cdot a_{<-1>}\ot ba_{<0>})\\
&=&\sum ba_{<0>}\cdot \gamma(c\cdot a_{<-1>})\\
&=& b\cdot (a\rhu \gamma)(c)\\
&=&b\cdot(\gamma\lhu a)(c)\\
&=&b\cdot(\gamma(c)\cdot a)\\
&=&\psi(c\ot b)\cdot a,
\end{eqnarray*}
and right $C^*$-linear since
\begin{eqnarray*}
\psi((c\ot b)\cdot c^*)&=&\psi(c\lhua c^*\ot b)\\
&=&b\cdot \gamma(c\lhua c^*)\\
&=& b\cdot (\gamma(c)\cdot c^*)\\
&=&\psi(c\ot b )\cdot c^*.
\end{eqnarray*}
Conversely, if $\psi\in V_5$, then $\gamma=g_4(\psi)\in V_4$. Indeed,
$\gamma$ is right $C^*$-linear, since
\begin{eqnarray*}
\gamma(c\lhua c^*)&=&\psi(c\lhua c^*\ot 1_A)\\
&=&\psi((c\ot 1_A)\cdot c^*)\\
&=& \psi(c\ot 1_A)\cdot c^*\\
&=&\gamma(c)\cdot c^*.
\end{eqnarray*}
For all $a\in A$, $c\in C$ we have
\begin{eqnarray*}
(a\rhu \gamma)(c)&=&\sum a_{<0>}\cdot \gamma(c\cdot a_{<-1>})\\
&=&\sum a_{<0>}\cdot \psi(c\cdot a_{<-1>}\ot 1_A)\\
(\psi {\rm~is~left~}A\hbox{-}{\rm linear})~~
&=&\sum \psi(c\cdot a_{<-1>}\ot a_{<0>})\\
&=&\psi((c\ot 1_A)\cdot a)\\
(\psi {\rm~is~right~}A\hbox{-}{\rm linear})~~
&=&\psi(c\ot 1_A)\cdot a\\
&=&(\gamma\lhu a)(c).
\end{eqnarray*}
proving that $\gamma$ is centralized by the left and right actions of
$A$. It is also clear that $f_4$ and $g_4$ are each others inverses,
and this finishes the proof of the fact that the $V_i$ are isomorphic
as $k$-modules. We leave it to the reader to prove that the $f_i$ and
$g_i$ send normalized elements to normalized elements.
\end{proof}

Combining \thref{2.1.3} and \prref{2.1.4}, we conclude that the (normalized)
elements of each of the five vector space $V_i$ can be used as tools to
construct a splitting of $\rho:\ 1_{\cal C}\to GF$, and, following the
philosophy of the introduction, then can be called integrals. We have reserved
this terminology for the elements $\gamma\in V_4$, because they are closely
related to Doi's total integrals (see \seref{3.1}).

\begin{definition}\delabel{2.1.4a}
A map $\gamma\in V_4$, that is, a $k$-linear map $\gamma:\ C\to \Hom(C,A)$
satisfying \eqref{2.1.1.4a} and \eqref{2.1.1.4e}
is called an $A$-integral for the
Doi-Hopf datum $(H,A,C)$. If $\gamma$ is normalized (that is, $\gamma$ satisfies
\eqref{2.1.1.4f}), then we call $\gamma$ a total $A$-integral.
\end{definition}

Our results can be summarized in the following Theorem, which is the
main result of this paper.

\begin{theorem}\thlabel{2.1.5}
{\bf (Maschke's Theorem for Doi-Hopf modules)}\\
Let $H$ be a bialgebra, and $(H,A,C)$ a Doi-Hopf datum.
If there exists a total $A$-integral
$\gamma:\ C\to \Hom(C,A)$ (or, equivalently, a normalized element in
$V_i$ ($i=1,\ldots,5$)), then the forgetful functor
$F:\ {}^C{\cal M}_A\to {\cal M}_A$ is separable. The converse holds if
$H$ is a Hopf algebra.\\
In this situation, the following assertions hold:\\
1) If a morphism $u:\ M\to N$ in ${}^C{\cal M}_A$ has a retraction
(resp. a section) in ${\cal M}_A$, then it has a retraction (resp. a section) in
${}^C{\cal M}_A$.\\
2) If $M\in {}^C{\cal M}_A$ is semisimple (resp. projective, injective)
as a right $A$-module, then
$M$ is also semisimple (resp. projective, injective) as an object in
${}^C{\cal M}_A$.
\end{theorem}

As we have remarked, the natural transformations not only form a $k$-module,
but even a $k$-algebra. The isomorphisms $f$ and $f_i$ then define an
algebra structure on each of the vector spaces $V_i$. The multiplication
can be described in a natural way for $i=3$ and $i=4$, and this is
what we will do next.\\
On $C^{\rm cop}\ot C$, we define an $H$-coaction
$$(c\ot d)\cdot h=c\ot d\cdot h$$
and this coaction makes $C^{\rm cop}\ot C$ into a right $H$-module coalgebra,
so that we can consider the $k$-algebra $\#(C^{\rm cop}\ot C,A)$
(see \eqref{Ko}).

\begin{lemma}\lelabel{2.1.6}
$V_3$ is a $k$-subalgebra of $\#(C^{\rm cop}\ot C,A)$. Normalized
elements in $V_3$ are right units of $V_3$, and consequently they are
idempotents.
\end{lemma}

\begin{proof}
Applying \eqref{2.1.1.3b}, we find the following formula for the
multiplication of $\theta,~\theta'\in V_3$:
\begin{eqnarray}
(\theta\bullet\theta')(c\ot d)&=&
\sum \theta(c_{(2)}\ot d_{(1)})_{<0>}
\theta'\bigl(c_{(1)}\ot d_{(2)}\cdot \theta(c_{(2)}\ot d_{(1)})_{<-1>}\bigr)
\nonumber\\
&=& \sum \theta(c_{(3)}\ot d)\theta'(c_{(1)}\ot c_{(2)})\eqlabel{2.1.6.1}
\end{eqnarray}
Using \eqref{2.1.6.1}, we can check easily that $\theta\bullet\theta'$
satisfies \eqref{2.1.1.3a} and \eqref{2.1.1.3b}. Indeed, for all
$c,d\in C$ and $a\in A$, we have
\begin{eqnarray*}
((\theta\bullet\theta')(c\ot d))\cdot a&=&
\sum \theta(c_{(3)}\ot d)a_{<0>}
\theta'(c_{(1)}\cdot a_{<-2>}\ot c_{(2)}\cdot a_{<-1>})\\
&=& \sum a_{<0>}\theta(c_{(3)}\cdot a_{<-2>}\ot d\cdot a_{<-1>})
\theta'(c_{(1)}\cdot a_{<-4>}\ot c_{(2)}\cdot a_{<-3>})\\
&=& \sum a_{<0>}(\theta\bullet\theta')(c\cdot a_{<-2>}\ot d\cdot a_{<-1>})
\end{eqnarray*}
and
\begin{eqnarray*}
\lefteqn{\sum d_{(2)}
\bigl((\theta\bullet\theta')(c\ot d_{(1)})\bigr)_{<-1>}\ot
\bigl((\theta\bullet\theta')(c\ot d_{(1)})\bigr)_{<0>}}\\
&=&
\sum d_{(2)}(\theta(c_{(3)}\ot d_{(1)}))_{<-1>}(\theta'(c_{(1)}\ot
c_{(2)}))_{<-1>}
\ot (\theta(c_{(3)}\ot d_{(1)}))_{<0>}(\theta'(c_{(1)}\ot c_{(2)}))_{<0>}\\
\eqref{2.1.1.3b}~~&=&
\sum c_{(3)}(\theta'(c_{(1)}\ot c_{(2)}))_{<-1>}\ot
\theta(c_{(4)}\ot d)(\theta'(c_{(1)}\ot c_{(2)}))_{<0>}\\
\eqref{2.1.1.3b}~~&=&
\sum c_{(1)} \ot \theta(c_{(4)}\ot d)\theta'(c_{(2)}\ot c_{(3)})
\end{eqnarray*}
>From \eqref{2.1.6.1}, it follows immediately that normalized elements in $V_3$
are right units.
\end{proof}

\begin{lemma}\lelabel{2.1.7}
$V_4$ is a $k$-subalgebra of $\Hom(C^{\rm cop},\#(C,A))$,
and $f_3:\ V_3\to V_4$ is a $k$-algebra isomorphism.
\end{lemma}

\begin{proof}
Left as an exercise to the reader.
\end{proof}

Now write $h=g\circ g_1\circ g_2:\ V_3\to V$.

\begin{theorem}\thlabel{2.1.8}
$h$ is a $k$-algebra homomorphism. If $H$ is a Hopf algebra, then $V$, $V_3$
and $V_4$ are isomorphic as $k$-algebras.
\end{theorem}

\begin{proof}
We have to show that $h$ is multiplicative. From the definitions of
$g,~ g_1$ and $g_2$, we easily find that $h(\theta)=\nu$ with
\begin{equation}\eqlabel{2.1.8.1}
\nu_M(c\ot m)=\sum m_{<0>}\cdot\theta(c\ot m_{<-1>})
\end{equation}
for all $M\in {\cal C}$. Now take $\theta'\in V_3$ and write
$h(\theta')=\nu'$. From \eqref{2.1.6.1}, we find that
$$(h(\theta\bullet\theta'))_M(c\ot m)=
\sum m_{<0>}\theta(c_{(3)}\ot m_{<-1>})\theta'(c_{(1)}\ot c_{(2)})$$
By definition $\nu\bullet\nu'=\nu'\circ \rho\circ \nu$, and
\begin{eqnarray*}
(\nu\bullet\nu')_M(c\ot m)&=&
\nu'_M(\rho_M(\nu_M(c\ot m)))\\
&=& \nu'_M(\sum c_{(1)}\ot \nu_M(c_{(2)}\ot m))\\
\eqref{2.1.8.1}~~~&=& \sum \nu_M(c_{(3)}\ot m_{<0>})\cdot
\theta'(c_{(1)}\ot c_{(2)})\\
\eqref{2.1.8.1}~~~&=& \sum m_{<0>}\theta(c_{(3)}\ot m_{<-1>})
\theta'(c_{(1)}\ot c_{(2)})
\end{eqnarray*}
and it follows that $h$ is multiplicative.
\end{proof}

\begin{example}\exlabel{2.1.9}\rm
Take $H=A=k$. The $k$-algebra $V_3$ then consists of maps
$\theta\in (C\ot C)^*$ satisfying
$$\sum c_{(1)}\ot \theta(c_{(2)}\ot d)=
\sum d_{(2)}\ot \theta(c\ot d_{(1)})$$
for all $c,d\in C$. The multiplication on $V_3$ is given by convolution,
and the normalized elements in $V_3$ are nothing else then the
coseparability idempotents in the sense of Larson \cite{La}. Thus
$C$ is a coseparable $k$-coalgebra if and only if the forgetful functor
${}^C{\cal M}\to k$-mod is separable.\\
If $k$ is a field, then this property is also equivalent to $C\ot L$ being
a cosemisimple coalgebra over $L$ for any field extension $L$ of $k$
(see \cite[Theorem 4.3]{CMZ3}). This is in fact the dual of the well-known
result that a $k$-algebra $A$ is separable over $k$ if and only if
$A\ot L$ is a semisimple $L$-algebra for any field extension $L$ of $k$.\\
We also remark that a method to construct coseparable coalgebras was
recently developed in \cite{CIMS}. Using an FRT type Theorem, one can
construct coseparability idempotents from solutions of the so-called
separability equation
$$R^{12}R^{23}=R^{23}R^{13}=R^{13}R^{12}$$
where $R\in \End(M\ot M)$, with $M$ a finite dimensional vector space.
\end{example}

\begin{remark}\relabel{2.1.10}\rm
Starting with a normalized element $\nu\in V_1$, we can construct a
splitting $\nu$ of $\rho$, and this $\nu$ is given by \eqref{2.1.3.1a}.
Of course one can also construct $\nu$ starting from
a normalized element in any of of the other four $V_i$. Let us do this
explicitely for $i=4$: take a total $A$-integral $\gamma\in V_4$. Thus
$\gamma:\ C\to \#(C,A)$, and the relation between $\gamma$ and the
corresponding $\nu\in V_1$ is the following:
$$\gamma(c)(d)=(\varepsilon_C\ot I_A)\nu(c\ot d\ot 1_A)$$
and consequently \eqref{2.1.3.1a} can be rewritten as
\begin{equation}\eqlabel{2.1.10.1}
\nu_M(c\ot m)=\sum m_{<0>}\gamma(c)(m_{<-1>})=m\trl \gamma(c)
\end{equation}
(see \eqref{trl} for the definition of the right $\#(C,A)$-action
on $M$). From this formula, we easily deduce the following result.
\end{remark}

\begin{proposition}\prlabel{2.1.11}
Let $H$ be a bialgebra.
For a normalized $A$-integral $\gamma:\ C\to \#(C,A)$, the following
statements are equivalent:\\
1) $\gamma$ is normalized;\\
2) for all $c\in C$, $a\in A$, we have
\begin{equation}\eqlabel{NO1}\eqlabel{2.1.11.1}
\sum (c_{(2)}\ot a)\trl \gamma(c_{(1)})=c\ot a
\end{equation}
3) for all $M\in {\cal C}$ and $m\in M$, we have
\begin{equation}\eqlabel{2.1.11.2}
\sum m_{<0>}\trl \gamma(m_{<-1>})=m
\end{equation}
\end{proposition}

\begin{proof}
$1)\Rightarrow 3)$ From \eqref{2.1.10.1} and the fact that $\nu$ splits $\rho$,
we find immediately that
$$m=\nu_M(\rho_M(m))=\sum m_{<0>}\trl \gamma(m_{<-1>})$$
$3)\Rightarrow 2)$ follows after we apply $3)$ with $M=C\ot A$.\\
$2)\Rightarrow 1)$ taking $a=1_A$, we find
$$
\sum c_{(3)}\cdot \gamma(c_{(1)})(c_{(2)})_{<-1>}\ot
\gamma(c_{(1)})(c_{(2)})_{<0>}=c\ot 1_A
$$
and \eqref{2.1.1.4f} follows after we apply $\varepsilon\ot I$.
\end{proof}

\subsection{The induction functor}\selabel{2.2}
In \cite{CMZ2}, the authors called
an element $z=\sum_i c_i\ot a_i\in C\ot A$ an integral for the Doi-Hopf datum
$(H,A,C)$ if
$$a z=z a$$
for all $a\in A$, or, equivalently, the map $\theta_A:\ A\to C\ot A$
given by $\theta_A(a)=za$ is left and right $A$-linear. From now on,
we wil call such a $z$ a {\it dual $A$-integral} for $(H,A,C)$. A
dual $A$-integral is called normalized if $\sum \varepsilon(c_i)a_i=1_A$.
It was shown in \cite{CMZ2} that the existence of dual $A$-integrals
is connected to the fact that the forgetful functor $F$ is Frobenius,
which means that $G$ is at once a left and right adjoint of $F$.
Now we will show that the existence of a normalized dual $A$-integral
is equivalent to the separability of the functor $G$.

\begin{theorem}\thlabel{2.2.1}
Let $(H,A,C)$ be a Doi-Hopf datum, and let $W$ and $W_1$ be the
$k$-modules consisting of respectively all natural transformations
$\theta:\ 1_{{\cal M}_A}\rightarrow FG $, and all dual $A$-integrals.
Then we have an isomorphism $f:\ W\to W_1$. The counit $\delta$ of
the adjoint pair $(F,G)$ is cosplit by $\delta\in W$ if and only if
the corresponding $A$-integral $f(\delta)$ is normalized.
Consequently the induction functor $G:\ {\cal M}_A\to {}^C{\cal M}_A$
is separable if and only if there exists a normalized dual integral
$z\in C\ot A$.
\end{theorem}

\begin{proof}
We define $f:\ W\to W_1$ as follows: $f(\theta)=\theta_A(1_A)$. We have
to show that $z=\theta_A(1_A)=\sum_i c_i\ot a_i\in W_1$.
>From the fact that $\theta_A$ is right $A$-linear, it follows that
$$\theta_A(a)=\sum_i c_i\cdot a_{<-1>}\ot a_ia_{<0>}$$
for all $a\in A$.
Let $N$ be an arbitrary right $A$-module. For $n\in N$, we define
$\phi_n:\ A\to N$ by $\phi_n(a)=na$, for all $a\in A$.
Using the naturality of $\theta$, we obtain the following commutative diagram
$$\begin{diagram}
A              &   \rTo^{\theta_A}   &  C\ot A                  \\
\dTo^{\phi_n}  &                     &  \dTo_{I_C \ot \phi_n}   \\
N              &  \rTo^{\theta_N}    &  C\ot N
\end{diagram}$$
Hence $\theta_N\circ \phi_n=(I_C\ot \phi_n)\circ \theta_A$.
Applying the diagram to $a\in A$ we find
\begin{eqnarray*}
\theta_N(na)&=& \sum c_{i}\cdot a_{<-1>}\ot \phi_n(a_ia_{<0>})\\
&=& \sum c_{i}\cdot a_{<-1>}\ot na_ia_{<0>}\\
&=& (\sum c_i\ot na_i)\cdot a
\end{eqnarray*}
and it follows that
\begin{equation}\eqlabel{2.2.1.1}
\theta_N (n)=\sum c_i\ot na_i
\end{equation}
Taking $N=A$ and $n=a$ we find
$$\sum c_i\ot aa_i=\sum c_i\cdot a_{<-1>}\ot a_ia_{<0>}$$
and this shows that $z$ is a dual $A$-integral.\\
Now define $g:\ W_1\to W$ by $g(\sum_i c_i\ot a_i)=\theta$,
with
$$\theta_N(n)=\sum_i c_i\ot na_i$$
for all $N\in {\cal M}_A$, and $n\in N$. Using the fact that
$z=\sum_i c_i\ot a_i$ is a dual $A$-integral, it follows easily that
$\theta_N$ is right $A$-linear. The naturality of $\theta$ can also be
verified easily: for any right $A$-linear map $\alpha:\ N\to N'$, we have
\begin{eqnarray*}
(I_C\ot \alpha)(\theta_N(a))&=&
\sum_i c_i\ot \alpha(na_i)\\
\sum_i c_i\ot \alpha(n)a_i&=&\theta_{N'}(\alpha(n))
\end{eqnarray*}
>From the definition of $f$ and $g$, it follows immediately that
$f\circ g= I_{W_1}$, and \eqref{2.2.1.1} implies $g\circ f=I_W$.\\
The statement about normalized elements is easy to verify, and the rest
of the Theorem then follows immediately from \thref{1.1}.
\end{proof}

\begin{remark}\relabel{2.2.2}\rm
$W$ is a $k$-algebra. The multiplication is given by the formula
$$\theta\circ \theta'=\theta'\circ \delta\circ \theta$$
The corresponding multiplication on $W_1$ is then the following:
$$\bigl(\sum_i c_i\ot a_i\bigr)\bullet \bigl(\sum_j c'_j\ot a'_j\bigr)=
\sum_{i,j} c'_j\ot\varepsilon(c_i)a_ia'_j$$
Obviously normalized dual $A$-integrals are left units, and idempotents
for this multiplication.
\end{remark}

\begin{example}\exlabel{2.2.3}\rm
We return to situation considered in \exref{2.1.9}.
>From \thref{2.2.1}, it follows immediately that the induction functor
$C\ot\bullet:\ {\cal M}_k\to {}^C{\cal M}$ is separable if and only if
there exists a $c\in C$ such that $\varepsilon_C(c)=1$. In particular,
if $C$ has a group-like element, then the functor $C\ot\bullet$
is separable, and this was pointed out already in \cite{Mi}.\\
If $C$ is a nonzero coalgebra over a field $k$, then $C\ot\bullet$ is
separable: take $c\neq 0\in C$, then $c=\sum\varepsilon(c_{(1)})c_{(2)}$,
and there exists a $d\in C$ with $\varepsilon(d)\neq 0$.
\end{example}

\subsection{Bialgebras and Relative separability}\selabel{2.3}
In \thref{2.1.3}, we have seen that the $k$-modules $V_1,\ldots,V_5$
are isomorphic to $V$ if $H$ is a Hopf algebra. Now assume that $H$ is
a bialgebra, not necessarily with an antipode. In this Section, we will
give a description of the $V_i$ in terms of natural transformations.
To this end, we first introduce a relative version of separability.\\
Let $(F,G)$ be a pair of adjoint functors between two categories
${\cal C}$ and ${\cal D}$. As usual, we write
$\rho:\ 1_{\cal C}\to GF$ for the unit of the adjunction. Let ${\cal B}$
be a third category such that we have functors $T:\ {\cal B}\to {\cal D}$
and $P:\ {\cal C}\to {\cal B}$, and a natural transformation
$\chi:\ TF\rightarrow P$:
$$
\begin{diagram}
{\cal C}  &         &\pile{ \rTo^F \\ \lTo_G} &       &{\cal D}  \\
    &\SE_P    &                         &\SW_T  &    \\
    &         &       {\cal B}                &       &    \\
\end{diagram}
$$
We call the functor $F$ {\it separable relative to} $({\cal B}, P,T,\chi)$
if there exists a natural transformation $\phi:\ PG\to T$ such that
\begin{equation}\eqlabel{2.3.1}
\chi\circ \phi(F)\circ P(\rho)=1_P
\end{equation}
the identity natural transformation on $P$.\\
If the functor $P$ is fully faithful, then relative separability implies
separability. Indeed, separability is equivalent to the splitting of the
unit $\rho$, and, because of the fact that $P$ is fully faithful, this
is equivalent to the fact that $P(\rho):\ P\to PGF$ splits, and this
happens if $P$ is separable relative to $({\cal B}, P,T,\chi)$.\\
Now we consider the following situation: ${\cal C}= {}^C{\cal M}_A$,
${\cal D}={\cal M}_A$, $F$ the functor forgetting the coaction,
${\cal B}={\cal M}_{A\# C^*}$, and the functors $P$ and $T$ defined as
in \seref{1.1}. $P$ is fully faithful if $C$ is projective as a $k$-module,
and we have the following diagram:
$$
\begin{diagram}
{}^C{\cal M}(H)_A  &         &\pile{ \rTo^{F} \\ \lTo_{G}} &      &{\cal
M}_A  \\
             &\SE_P    &                                     &\SW_T &      \\
             &         &    {\cal M}_{A\# C^*}                     &      &
\\
\end{diagram}
$$
We define $\chi:\ TG\to P$ as follows:
$\chi_M:\ M\ot_A \#(C,A)\to M$ is given by
$$\chi_M(m\ot f)=m\trl f$$
for all $M\in {\cal M}_A$, $m\in M$ and $f\in\Hom(C,A)$.

\begin{theorem}\thlabel{2.3.1}
Let $H$ be a bialgebra, and $(H,A,C)$ a Doi-Hopf datum.
With notations as above, the $k$-module $V'$ consisting of all
natural transformations $\phi:\ PG\rightarrow T$ is isomorphic to $V'_5$,
the $k$-module consisting of all $(A,A\#C^*)$-bimodule maps
$\psi:\ C\ot A\to \Hom(C,A)$. If $C$ is projective as a $k$-module, then
$V'\cong V'_5$ is isomorphic to each of the $V_i$, and, consequently,
$F$ is separable relative to $({\cal B},P,T,\chi)$ if and only if there exists a
normalized integral in one of the $V_i$'s.
\end{theorem}

\begin{proof}
We will construct an isomorphism $g':\ V'_5\to V'$.
Let $\psi:\ C\ot A\to \Hom(C,A)$ be an integral in $V_5$. We define
$$\phi=g'(\psi):\ PG\rightarrow T$$
by
$$\phi_M:\ C\ot M\to M\ot_A\#(C,A),~~~
\phi_M(c\ot m)=m\ot\psi(c\ot 1_A)$$
for all $M\in {\cal M}_A$, $c\in C$ and $m\in M$.
We will prove that $\phi$ is a natural transformation from $PG$
to $T$. First, $\phi_M$ is right $A\# C^*$-linear. For all
$a\in A$, $c\in C$, $c^*\in C^*$ and $m\in M$ we have
\begin{eqnarray*}
\phi_M((c\ot m)\cdot a)&=&\sum \phi_M(c\cdot a_{<-1>}\ot m\cdot a_{<0>})\\
&=&\sum m\cdot a_{<0>}\ot_A \psi(c\cdot a_{<-1>}\ot 1_A)\\
&=& \sum m\ot_A a_{<0>}\cdot \psi(c\cdot a_{<-1>}\ot 1_A)\\
&=&m\ot_A  (a\cdot\psi)(c\ot 1_A)\\
&=&m\ot_A \psi(c\ot a)\\
&=&m\ot_A \psi(c\ot 1_A)\cdot a\\
&=&\phi_M(c\ot m)\cdot a
\end{eqnarray*}
and
\begin{eqnarray*}
\phi_M((c\ot m)\cdot c^*)&=&\phi_M(c\lhua c^*\ot m)\\
&=&m\ot_A \psi(c\lhua c^*\ot 1_A)\\
&=&m\ot \psi(c\ot 1_A)\cdot c^*\\
&=&\phi_M(c\ot m)\cdot c^*
\end{eqnarray*}
The naturality condition can be verified as follows.
Take $M,N\in{\cal M}_A$ and a right $A$-linear map $f:M\to N$. Then
\begin{eqnarray*}
(f\ot_AI_{\#(C,A)})\circ \phi_M(c\ot m)&=&
(f\ot_AI_{\#(C,A)})(m\ot_A \psi(c\ot 1_A))\\
&=& f(m)\ot_A\psi(c\ot 1_A)\\
&=&\phi_N(c\ot f(m))\\
&=&\phi_N\circ (I_C\ot f).
\end{eqnarray*}
$f':\ V'\to V'_5$ is defined as follows: for any natural transformation
$\phi:\ PG\rightarrow T$, we take
$$f'(\phi)=\psi=\phi_A:\ C\ot A\to A\ot_A\#(C,A)\cong \#(C,A)$$
We have to show that $\psi$ is an integral in $V'_5$.\\
$\psi$ is right $A\# C^*$-linear by definition.
The next thing to show is that $\psi$ is left $A$-linear. Let
$N\in{\cal M}_A$ and
$n\in N$, and consider the map
$$\varphi_n:\ A\to N;~~\varphi_n(a)=na$$
for all $a\in A$. It is clear that $\varphi_n$ is right $A$-linear, and,
by the naturality of of $\phi$,
we have the following commutative diagram
$$
\begin{diagram}
C\ot A                    & \rTo^{\phi_A} & A\ot_A\#(C,A)                   \\
\dTo^{I_C\ot\varphi_n}  &               & \dTo_{\varphi_n\ot I_{\#(C,A)}} \\
C\ot N                    & \rTo^{\phi_N} & N\ot_A\#(C,A)                   \\
\end{diagram}
$$
Applying the diagram to $c\ot a\in C\ot A$, we obtain
\begin{equation}\eqlabel{fu}
(\varphi_n\ot I_{\#(C,A)})(\phi_A(c\ot a))=\phi_N(c\ot na),
\end{equation}
for all $c\in C$, $a\in A$. Taking $a=1$, we find
\begin{equation}\eqlabel{gen}
\phi_N(c\ot n)=n\trl \psi(c\ot 1).
\end{equation}
If we apply \eqref{fu} for $N=A$ and $n=b\in A$ we get
$$b\cdot\psi(c\ot a)=\psi(c\ot ba),$$
and this proves that $\psi$ is left $A$-linear.\\
It is easy to see that $f'$ and $g'$ are each others inverses.
Let $\phi=g'(\psi)$. Then $\phi_A(c\ot a)=a\psi(c\ot 1_A)=\psi(c\ot a)$,
so $\phi_A=f'(g'(\psi))=\psi$, and $f'\circ g'=I_{V_3}$.\\
Conversely, take $\phi\in V'$, and write $\tilde{\phi}=(g'\circ f')(\phi)$.

Then
\begin{eqnarray*}
\tilde{\phi}_N(c\ot n)&=& n\ot \phi_A(c\ot 1_A)\\
&=& (\varphi_n\ot I_{\#(C,A)})(\phi_A(c\ot 1_A))\\
\eqref{fu}~~~&=& \phi_N(c\ot n)
\end{eqnarray*}
and this proves that $g'\circ f'=I_{V'}$.\\
Finally $\chi\circ\phi(F)\circ P(\rho)=1_P$ if and only if
$$\chi_M\circ \phi_{C\ot M}\circ\rho_M=I_M$$
for all $M\in {}^C{\cal M}_A$, and this condition in fact means that
$$\sum m_{<0>}\trl \psi(m_{<-1>}\ot 1_A)=m$$
for all $m\in M$, and this equivalent to $M$ being normalized, by
\prref{2.1.11}.
\end{proof}

\section{Applications and examples}\selabel{3}
In this Section, we will make clear how \thref{2.1.5} generalizes some
existing versions of Maschke's Theorem. We will also give some new
applications and examples.

\subsection{Previous types of integrals for Doi-Hopf modules}\selabel{3.2}
Let $(H,A,C)$ be a Doi-Hopf datum, where $H$ is a
Hopf algebra with a bijective antipode $S$,
and $C$ is projective as a $k$-module.
According to \cite{CMZ3}, an integral is an $(A,A\# C^*)$-linear map
$$\psi:\ C\ot A \to C^*\ot A$$
where the $(A,A\# C^*)$-bimodule structure on $C^*\ot A$ is given by
the following formulas:
\begin{eqnarray*}
b\cdot (c^*\ot a)&=&c^*\ot ba\\
(c^*\ot a)\cdot b&=&\sum \lan c^*,?\cdot S^{-1}(b_{<-1>})\ran\ot ab_{<0>}\\
(c^*\ot a)\cdot d^*&=&c^**d^*\ot a
\end{eqnarray*}
for all $c^*,~d^*\in C^*$ and $a,~b\in A$. Such an integral is called
normalized if
$$\sum \left( (c_{(2)}\ot a)\cdot (c_{(1)})^{[0]}\right)\cdot
(c_{(1)})^{[-1]}=c\ot a$$
for all $c\in C$, $a\in A$ (we denoted
$\psi(c\ot 1_A)=\sum c^{[-1]}\ot c^{[0]}$\ ).\\
If $\psi:\ C\ot A\to C^*\ot A$ is a normalized integral in the sense of
\cite{CMZ3} then $\gamma:\ C\to \Hom(C,A)$, defined by
$\gamma(c)=i\circ f\circ \psi(c\ot 1_A)$, for all $c\in C$, is a normalized
$A$-integral in our sense. Here $f:C^*\ot A\to A\# C^*$ is
the $(A,A\# C^*)$-bimodule map given by the formula
$$f(c^*\ot a)=\sum a_{<0>}\# a_{<-1>}\cdot c^*$$
A straightforward computation shows that $f$ is $(A, A\# C^*)$-linear.
Let us check that $\gamma$ is
normalized. For all
$c\in C$ and $a\in A$, we have that
\begin{eqnarray*}
\sum (c_{(2)}\ot a)\trl \gamma(c_{(1)}) &=&
\sum c_{(3)}\cdot \gamma(c_{(1)})(c_{(2)})_{<-1>}\ot
a\gamma(c_{(1)})(c_{(2)})_{<0>}\\
&=& \sum \left( (c_{(2)}\ot a)\cdot (c_{(1)})^{[0]}\right)\cdot
(c_{(1)})^{[-1]}\\
&=&c\ot a.
\end{eqnarray*}
and it follows from \eqref{NO1} that $\gamma$ is normalized.\\
Thus, all the examples of normalized integrals given in \cite{CMZ3}
provide examples of normalized $A$-integrals in our new sense.
For example, let $(k,k,C)$ be the Doi-Hopf datum  where $C={\rm M}^n(k)$,
is the
$n\times n$ matrix coalgebra, that is, $C$ is the dual of the $n\times n$
matrix algebra
${\rm M}_n(K)$. Let
$\{e_{ij},\varepsilon_{ij}|1\leq i,j\leq n\}$ be the canonical dual basis
for $C$, and let
$\mu=(\mu_{ij})$ be an arbitrary $n\times n$-matrix. Then the map
$\gamma_{\mu}:\ C\to C^{*}$ defined by
$$\gamma_{\mu}(\varepsilon_{ij})= \sum_{k=1}^n\mu_{kj}e_{ki}$$
is a $k$-integral map for $(k,k,C)$.
Furthermore, $\gamma_{\mu}$ is normalized if and only if Tr$(\mu)=1$ (see
Example 2.3 of
\cite{CMZ3}).

\subsection{Relative Hopf modules}\selabel{3.1}
\subsubsection{Doi's total integrals and the forgetful functor}
Let $H$ be a projective Hopf algebra and $A$ a left $H$-comodule algebra.
$H$ is
an $H$-module coalgebra, and the category
${}^H{\cal M}(H)_A$ is nothing else then the category of relative Hopf modules
${}^H{\cal M}_A$. We recall that a total integral
$\varphi:\ H\to A$ is left $H$-colinear map such that
$\varphi(1_H)=1_A$ (see \cite{D0}).\\
Several results about the separability of the forgetful functor
$F:\ {}^H{\cal M}_A\to {\cal M}_A$ have appeared in the literature.
Doi \cite[Theorem 1.7]{D0} proved that the following condition is sufficient:\\
$\bullet$ $H$ is commutative and $\varphi:\ H\to A$ is a total integral
such that
$\varphi(H)\subset Z(A)$.\\
In \cite[Theorem 1]{D1} the following two sufficient conditions for $F$ to be
separable are given:\\
$\bullet$ $H$ is involutory (i.e. $S^2=I_H$), $\varphi(H)\subseteq Z(A)$ and
$\varphi(hk)=\varphi(kh)$ for all $h,k\in H$;\\
$\bullet$ $A$ is faithful as a $k$-module and $\varphi(H)\subseteq k$.\\
and in this paper Doi asks if we can prove a Maschke type Theorem
(in our language: a separable functor Theorem) for relative Hopf modules
under more general assumptions. \thref{2.1.5} gives an answer this
problem:
\begin{corollary}\colabel{3.1.1a}
Let $H$ be a projective $k$-Hopf algebra, and $A$ a left
$H$-comodule algebra. Then the following statements are equivalent:\\
1) the forgetful functor $F:\ {}^H{\cal M}_A\to {\cal M}_A$ is separable;\\
2) there exists a normalized $A$-integral $\gamma:\ H\to \Hom(H, A)$.
\end{corollary}
We will now investigate the relation between Doi's total integrals and our
total $A$-integrals. This will explain our terminology, and we will also
prove that
the forgetful functor is separable if and only if there exists a total integral
$\varphi$ such that the image of $\rho\circ\varphi$ is contained in the center
of $H\ot A$. We will see that this last condition is implied by Doi's conditions
mentioned above (see Corolaries \ref{co:3.1.3} and \ref{co:3.1.4}).

\begin{proposition}\prlabel{3.1.1}
Let $H$ be a $k$-projective bialgebra and $A$ be a left $H$-comodule algebra.
If $\gamma:\ H\to \Hom(H,A) $ is a (normalized) $A$-integral for $(H,A,H)$ then
$$\varphi_{\gamma}:\ H\to A, \quad  \varphi_\gamma (h):=\gamma (h)(1_H),$$
for all $h\in H$, is a total integral in the sense of Doi \cite{D0}.
\end{proposition}

\begin{proof}
$\varphi_{\gamma}$ is right $H^*$-linear since
\begin{eqnarray*}
\varphi _\gamma (h\lhua h^*)&=& \gamma(h\lhua h^*)(1_H)\\
&=&(\gamma(h)\cdot h^*)(1_H)\\
&=&\sum \gamma(h)(1_H)_{<0>}\lan h^*, \gamma(h)(1_H)_{<-1>}\ran\\
&=& (\gamma(h)(1_H))\cdot h^*\\
&=& \varphi _\gamma (h) \cdot h^*
\end{eqnarray*}
for all $h\in H,~h^*\in H^*$. Furthermore, if $\gamma$ is normalized, then
we obtain from \eqref{2.1.1.4f}
$$\varphi_\gamma (1_H)=\gamma(1_H)(1_H)=\varepsilon(1_H)1_A=1_A.$$
proving that $\varphi_{\gamma}$ is a total integral.
\end{proof}

Conversely, let $\varphi:\ H\to A$ be a (total) integral, and define
\begin{equation}\eqlabel{3.1.2.1}
\gamma ^\varphi:\ H\to \Hom(H,A),~~~\gamma ^\varphi (h)(k)=\varphi(hS(k)),
\end{equation}
for all $h$, $k\in H$. We will now present some necessary and sufficient
condition for $\gamma ^\varphi$ to be a (normalized) $A$-integral in $V_4$.

\begin{theorem}\label{total}\thlabel{3.1.2}
Let $A$ be a left $H$-comodule algebra and $\varphi:\ H \to A $ an integral.
If
\begin{equation}\eqlabel{3.1.2.2}
\rho \Bigl (\varphi (H) \Bigl) \subseteq Z(H\otimes A),
\end{equation}
the center of the tensor product of $H$ and $A$, then the map $\gamma^\varphi$
defined by \eqref{3.1.2.1} is an
$A$-integral for $(H,A,H)$.\\
Conversely, if $H$ is projective as a $k$-module, and $\gamma^\varphi$ is an
$A$-integral, then \eqref{3.1.2.2} holds.\\
Finally, $\varphi$ is a total integral if and only if $\gamma^\varphi$ is
normalized $A$-integral.
\end{theorem}

\begin{proof}
First we remark that \eqref{3.1.2.2} is equivalent to
\begin{equation}\eqlabel{3.1.2.3}
(g\otimes 1_A)(\rho\varphi(h)) = (\rho\varphi(h)) (g\otimes 1_A)
\end{equation}
for all $g,h\in H$ and
\begin{equation}\eqlabel{3.1.2.4}
\varphi(h)\subset Z(A)
\end{equation}
Now assume that \eqref{3.1.2.2} holds.
Then for any $h,g\in H,\ h^*\in H^*$ we have
\begin{eqnarray*}
(\gamma^\varphi(h)\cdot h^*)(g) &=& \sum
\gamma^\varphi(h)(g_{(1)})_{<0>}{}\lan h^*,g_{(2)}
\gamma^\varphi(h)(g_{(1)})_{<-1>}\ran\\
&=& \sum \varphi(hS(g_{(1)}))_{<0>}{}\lan
h^*,g_{(2)}\varphi(hS(g_{(1)}))_{<-1>}\ran\\
\mbox{(by assumption)}\qquad &=& \sum
\varphi((hS(g_{(1)}))_{(2)}){}\lan h^*,g_{(2)}(hSg_{(1)})_{(1)}\ran\\
\mbox{($\varphi$ is $H$-colinear)}\qquad
&=& \sum \varphi((hS(g_{(1)}))_{(2)}){}\lan h^*,(hSg_{(1)})_{(1)}g_{(2)}\ran\\
&=& \sum \varphi(h_{(2)}S(g_{(1)})){}\lan h^*,h_{(1)}Sg_{(2)}g_{(3)}\ran\\
&=& \sum \varphi(h_{(2)}S(g)){}\lan h^*,h_{(1)}\ran\\
&=& \varphi((h\leftharpoonup h^*)S(g))\\
&=& \gamma^\varphi(h\leftharpoonup h^*)(g).
\end{eqnarray*}
So $(\gamma^\varphi(h)\cdot h^*)=\gamma^\varphi(h\leftharpoonup h^*)$,
for any $h\in H,\ h^*\in H^*$, and $\gamma^\varphi$ is right $H^*$-linear.\\
Now for any $a\in A,\ g,\ h\in H$,
\begin{eqnarray*}
(a\rhu \gamma^\varphi)(h)(g) &=& \sum (a_{<0>}\cdot
\gamma^\varphi(ha_{<-1>}))(g)\\
&=& \sum a_{<0>}[\gamma^\varphi(ha_{<-2>})(ga_{<-1>})]\\
&=& \sum a_{<0>}\varphi(ha_{<-2>}S(a_{<-1>})S(g))\\
&=& a\varphi(hS(g))\\
&=& a\gamma^\varphi(h)(g)\\
\mbox{(by assumption)}\qquad &=& \gamma^\varphi(h)(g)a\\
&=& (\gamma^\varphi\lhu a)(h)(g).
\end{eqnarray*}
and $a\rhu \gamma^\varphi = \gamma^\varphi\lhu a$ for any $a\in A$. This
proves that
$\gamma^\varphi$ is an $A$-integral.\\
Conversely, assume that $H$ is projective as a $k$-module and that
$\gamma^\varphi$ is an
$A$-integral. Then for any $h\in H,\ h^*\in H^*$, we have
$$\gamma^\varphi(h\leftharpoonup h^*) = \gamma^\varphi(h)\leftharpoonup h^*$$
Thus for any $g,\ h\in H,\ h^*\in H^*$,
\begin{eqnarray*}
\sum \varphi(h_{(2)}){}\lan h^*, h_{(1)}g\ran {} &=&
\sum \varphi(h_{(2)}g_{(2)}S(g_{(3)})){}\lan h^*, h_{(1)}g_{(1)}\ran\\
&=& \sum \varphi(((hg_{(1)})\leftharpoonup h^*)S(g_{(2)}))\\
 &=&\sum \gamma^\varphi((hg_{(1)})\leftharpoonup h^*)(g_{(2)})\\
\mbox{($\gamma^\varphi$ is integral)}\qquad
&=&\sum (\gamma^\varphi(hg_{(1)})\leftharpoonup h^*)(g_{(2)})\\
\mbox{(by definition of right $H^*$-module)}\qquad
&=&\sum [\gamma^\varphi(hg_{(1)})(gg_{(1)})]_{<0>} {}
\lan h^*, g_{(3)}[\gamma^\varphi(hg_{(1)})(g_{(2)})]_{<-1>}\ran\\
&=&\sum [\varphi(h)]_{<0>}{}\lan h^*, g[\varphi(h)]_{<-1>}\ran \\
&=&\sum \varphi(hg_{(1)}) {}\lan h^*, gh_{(1)}\ran .
\end{eqnarray*}
or
$$\sum \varphi(h_{(2)}){}\lan h^*, h_{(1)}g\ran {} =
\sum \varphi(hg_{(1)}) {}\lan h^*, gh_{(1)}\ran$$
for any $h,\ g\in H,\ h^*\in H^*$. Using the fact that $H$ is projective
as a $k$-module, we obtain
$$\sum h_{(1)}g\otimes \varphi(h_{(2)})=
\sum gh_{(1)}\otimes \varphi(h_{(2)})$$
and \eqref{3.1.2.3} follows from the fact that $\varphi$ is $H$-colinear.\\
For all $a\in A$, we have that
$a\rhu \gamma^\varphi = \gamma^\varphi \lhu a$. Therefore we find for any
$h\in H$ that
\begin{eqnarray*}
\varphi(h)a &=& \gamma^\varphi(h)(1)a\\
&=& (\gamma^\varphi(h)\cdot a)(1)\\
&=& (\gamma^\varphi \lhu a)(h)(1)\\
&=& (a\rhu \gamma^\varphi)(h)(1)\\
&=& \sum [a_{<0>}\gamma^\varphi(h a_{<-1>})](1)\\
&=& \sum a_{<0>}\gamma^\varphi(h a_{<-2>})(a_{<-1>})\\
&=& \sum a_{<0>}\varphi(h a_{<-2>}S(a_{<-1>}))\\
&=& a\varphi(h).
\end{eqnarray*}
So $\varphi(h)\in Z(A)$ for any $h\in H$, and \eqref{3.1.2.4} follows. As
we remarked
earlier, (\ref{eq:3.1.2.3}-\ref{eq:3.1.2.4}) imply \eqref{3.1.2.2}.\\
It is routine to check that for any left $H$-comodule map $\varphi:\ H\to A $
$$\varphi_{\gamma^{\varphi}} = \varphi$$
So if $\gamma^{\varphi}$ is a normalized integral then $\varphi$ is total by
\prref{3.1.1}. Finally, if $\varphi$ is a total integral, then
$$\gamma^\varphi(h_{(1)})(h_{(2)})=\varphi(h_{(1)}Sh_{(2)}) =
\varepsilon(h)\varphi(1) =
\varepsilon(h)1_A$$
for all $h\in H$, and $\gamma^\varphi$ is normalized.
\end{proof}

\begin{corollary}\colabel{3.1.3}
Let $\varphi:\ H\to A$ be left $H$-colinear. If
$$\varphi(gh)=\varphi(hS^2(g))\ \mbox{and}\ \varphi(H)\subseteq Z(A)$$
for all $g,\ h\in H$, then $\gamma^\varphi$ is an $A$-integral.
\end{corollary}

\begin{proof}
It follows immediately from our assumptions that
\begin{equation}\eqlabel{3.1.3.1}
\varphi(h)= \varphi(h1_H)=\varphi(1_HS^2h)=\varphi(S^2h)
\end{equation}
for all $h\in H$.
Applying $\rho$, and using the fact that $\varphi$ is $H$-colinear, we obtain
\begin{equation}\eqlabel{3.1.3.2}
\sum h_{(1)}\otimes \varphi(h_{(2)})=\sum S^2(h_{(1)})\otimes
\varphi(S^2(h_{(2)})) = \sum S^2(h_{(1)})\otimes \varphi(h_{(2)})
\end{equation}
and
\begin{eqnarray*}
\sum h_{(1)}S^2(g)\otimes \varphi(h_{(2)})&=&
\sum S^2(h_{(1)}g)\otimes \varphi(h_{(2)})\\
&=& \sum S^2(h_{(1)})S^2(g_{(1)})\otimes
\varphi(h_{(2)}g_{(2)}S(g_{(3)}))\\
&=& \sum S^2(h_{(1)})S^2(g_{(1)})\otimes
\varphi(S(g_{(3)})S^2(h_{(2)}g_{(2)}))\\
&=& \sum S^2(g_{(5)}) S(g_{(4)})S^2(h_{(1)}) S^2(g_{(1)})\otimes
\varphi(S(g_{(3)})S^2(h_{(2)})S^2(g_{(2)}))\\
&=&\sum S^2(g_{(3)})(S(g_{(2)})S^2(h)S^2(g_{(1)}))_{(1)}\otimes
\varphi((S(g_{(2)})S^2(h)S^2(g_{(1)}))_{(2)})\\
&=&\sum (S^2(g_{(3)})\otimes 1_A)\rho\varphi(S(g_{(2)})S^2(h)S^2(g_{(1)}))\\
&=&\sum (S^2(g_{(3)})\otimes 1_A)\rho\varphi(hg_{(1)}S(g_{(2)}))\\
&=&(S^2(g)\otimes 1_A)\rho\varphi(h)\\
&=&\sum S^2(g)h_{(1)}\otimes \varphi(h_{(2)}).
 \end{eqnarray*}
For all $h,\ g\in H$. It follows that
$$\sum h_{(1)}S^2(g)\otimes \varphi(h_{(2)})=\sum S^2(g)h_{(1)}\otimes
\varphi(h_{(2)})$$
and
\begin{eqnarray*}
\sum S(g)h_{(1)}\otimes \varphi(h_{(2)}) &= &
\sum S(g_{(3)})h_{(1)}S^2(g_{(2)})S(g_{(1)})\otimes
\varphi(h_{(2)})\\
&= & \sum S(g_{(3)})S^2(g_{(2)})h_{(1)}S(g_{(1)})\otimes
\varphi(h_{(2)})\\
&= & \sum h_{(1)}S(g)\otimes
\varphi(h_{(2)}),
\end{eqnarray*}
so
\begin{eqnarray*}
\sum gh_{(1)}\otimes \varphi(h_{(2)}) &= &
\sum g_{(1)}h_{(1)}S(g_{(2)})g_{(3)}\otimes\varphi(h_{(2)})\\
&= & \sum g_{(1)}S(g_{(2)})h_{(1)}g_{(3)}\otimes\varphi(h_{(2)})\\
&= & \sum h_{(1)}g\otimes\varphi(h_{(2)})
\end{eqnarray*}
This proves that $(g\otimes 1_A)\rho\varphi(h)=\rho\varphi(h)(g\otimes 1_A)$.
Now observe that $\varphi(H)\subseteq Z(A)$, and
$$\rho\Bigl (\varphi(H)\Bigl) \subseteq Z(H\otimes A)$$
and it follows from \thref{3.1.2} that $\gamma^\varphi$ is an integral.
\end{proof}

As a special case, we recover the following result of Doi \cite{D1}:

\begin{corollary}\colabel{3.1.4}
Let $H$ be a Hopf algebra, $A$ a left $H$-comodule algebra,
and $\varphi:\ H\to A$ be a (total) integral. Assume that one of
the following two conditions holds:\\
1) $H$ is involutory (i.e. $S^2=I_H$), $\varphi(H)\subseteq Z(A)$ and
$\varphi(hk)=\varphi(kh)$ for all $h,k\in H$;\\
2) The antipode of $H$ is bijective and $\varphi(H)\subseteq k$.\\
Then $\gamma^\varphi$ is a (normalized) $A$-integral, and the forgetful
functor ${}^H{\cal M}_A\to {\cal M}_A$ is separable.
\end{corollary}

\begin{proof}
1) follows immediately from \coref{3.1.3}.\\
2) Using the fact that $\varphi$ is left $H$-colinear, we find for all
$h\in H$ that
$$\sum h_{(1)}\otimes \varphi( h_{(2)})=1_H\otimes \varphi(h)$$
It is obvious that $\rho\varphi(H)\subseteq Z(H\otimes A)$, and the result
follows
again from \coref{3.1.3}.
\end{proof}

\subsubsection{Hopf Galois extensions and the induction functor}
Let $H$ be a Hopf algebra, and $A$ be a left $H$-comodule algebra. Recall
that the subalgebra of coinvariants is given by the formula
$$B=A^{{\rm co}H}=\{a\in A~|~\rho(a)=1\ot a\}$$
It is well-known that we have an adoint pair of functors (see e.g.
\cite{CR})
$$\bullet\ot_B A:\ {\cal M}_B\to {}^H{\cal M}_A~~;~~
(\bullet)^{{\rm co}H}:\ {}^H{\cal M}_A\to {\cal M}_B$$
Let $\delta$ be the counit of this adjunction. We know that
$H\ot A$ is a relative Hopf module, and, by definition, $A/B$ is a Hopf
Galois extension if the map $\delta_{H\ot A}:\ A\ot_B A\to H\ot A$
is an isomorphism. This adjunction map is given by the formula
$$\delta_{H\ot A}(a\ot b)=\sum a_{<-1>}\ot a_{<0>}b$$
Here we have to be careful with the structure maps on $H\ot A$. They are
not given by (\ref{eq:1.3}-\ref{eq:1.4}), but by
\begin{eqnarray}
(h\ot b)\cdot a&=& \sum h\ot ba\eqlabel{3.1.5.1}\\
\rho'(h\ot b)&=& \sum h_{(1)}b_{<-1>}\ot h_{(2)}\ot b_{<0>}\eqlabel{3.1.5.2}
\end{eqnarray}
It is well-known that the isomorphism
$$f:\ H\ot A\to H\ot A~~;~~f(h\ot a)=\sum ha_{<-1>}\ot a_{<0>}$$
translates the structures (\ref{eq:3.1.5.1}-\ref{eq:3.1.5.2}) into
(\ref{eq:1.3}-\ref{eq:1.4}), and consequently $A/B$
is a Hopf Galois extension if and only if the map
$$f\circ\delta_{H\ot A}=\delta'_{H\ot A}:\ A\ot_B A\to H\ot A~~;~~
\delta'_{H\ot A}(a\ot b)=\sum b_{<-1>}\ot ab_{<0>}$$
is an isomorphism. $\delta'_{H\ot A}$ is a map in ${}^H{\cal M}_A$ if we
give $H\ot A$ the usual structures (\ref{eq:1.3}-\ref{eq:1.4}). We now have
the following result:

\begin{theorem}\thlabel{3.1.5}
Let $H$ be a Hopf algebra, and $A$ a left $H$-comodule algebra, and
$B=A^{{\rm co}H}$. If $A$ is a separable extension of $B$, then
the induction functor $H\ot\bullet:\ {\cal M}_A\to {}^H{\cal M}_A$
is separable. The converse property holds if $A/B$ is a Hopf Galois
extension.
\end{theorem}

\begin{proof}
We have remarked already that $\delta'_{H\ot A}$ is right $A$-linear, and
it can be verified easily that it is also left $A$-linear. Let
$e=\sum e^1\ot e^2$ be a separability idempotent in $A\ot_B A$. Then
$e$ has the following properties:
$$\sum ae^1\ot e^2=\sum e^1\ot e^2a~~{\rm and}~~\sum e^1 e^2=1_A$$
for all $a\in A$.
Consider $z=\sum e^2_{<-1>}\ot e^1e^2_{<0>}=\delta'_{H\ot A}(e)\in H\ot A$.
>From the fact that $\delta'_{H\ot A}$ is left and right $A$-linear, it
follows immediately that $az=za$, so $z$ is a dual $A$-integral in the
sense of \seref{2.2}. Also the normalization condition follows easily,
since
$$\sum \varepsilon(e^2_{<-1>}) e^1e^2_{<0>}=\sum e^1 e^2=1_A$$
and the first statement follows from \thref{2.2.1}.\\
If $A/B$ is a Hopf Galois extension, then we proceed in a similar way,
but using the inverse of $\delta'_{H\ot A}$.
\end{proof}

\begin{remark}\rm
In the literature we often find the category ${\cal M}(H)^H_A$ instead of
${}^H{\cal M}_A$, where $A$ is now a right $H$-comodule algebra. Adapting
our results, we find that the functor $\bullet\ot H:\ {\cal M}_A\to
{\cal M}^H_A$ is separable if $A/A^{{\rm co}H}$ is separable, and the
converse holds if $A$ is a Hopf Galois extension of $A^{{\rm co}H}$.
\end{remark}

\subsection{Classical integrals}\selabel{3.3}
Let $H$ be a Hopf algebra, and assume that $H$ is flat as a
$k$-module. Recall from \cite{S} that $\varphi\in H^*$ is called
a left integral on $H^*$ if $h^**\varphi=\lan h^*,1_H\ran\varphi$ for all
$h^*\in H^*$,
or, equivalently, if $\varphi:\ H\to k$ is left (or right) $H$-colinear.\\
Now suppose that there exists a (total) $k$-integral $\gamma$ for
$(k,k,H)$, or, equivalently, a map $\theta:\ H\ot H\to k$ in $V_3$
(see \exref{2.1.9}). The map $\varphi=i(\gamma)\in H^*$ defined by
$$\lan\varphi,h\ran=\theta(h\ot 1_H)=\gamma(h)(1_H)$$
is a left integral. Conversely, if $\varphi\in\int_{H^*}^l$, the $k$-module
consisting of classical integrals on $H$, then
$p(\varphi)=\gamma:\ H\to H^*$ given by
$$\gamma(h)(g)=\varphi(hS(g))$$
is a $k$-integral. This can be proved directly, but it also follows from
\thref{3.1.2}. So we have maps
$$i:\ \int_{H^*}^l\to V_3\cong V_4~~{\rm and}~~
p:\ V_4\cong V_3\to \int_{H^*}^l$$
and it can be seen easily that $p$ is a left inverse of $i$.\\
Surprisingly, $i$ is not an isomorphism. To see this, let $H$ be a finite
dimensional Hopf algebra over a field $k$. It is well-known that
${\rm dim}_k(\int_{H^*}^l)=1$. On the other hand, a $k$-integral $\gamma$ is
nothing else then a right $H^*$-linear map $\gamma:\ H\to H^*$.
Now if $t\in H$ is a left integral, then $H=t\klop H^*$ (\cite{S}), hence
${\rm dim}_k(V_4)={\rm dim}_k(\Hom_{H^*}(H,H^*))=n$, and this shows that
the space of $k$-integrals $V_4$ is larger then the classical $\int_{H^*}^l$.\\
$\im(i)$ is the subspace of $V_3$ consisting of all $\theta\in V_3$
satisfying $(i\circ p)(\theta)=\theta$, or
$$\theta(h\ot k)=\theta(hS(k)\ot 1_H)$$
for all $h,k\in H$. This condition is equivalent to
\begin{equation}\eqlabel{3.3.1.1}
\sum\theta(hl_{(1)}\ot kl_{(2)})=\theta(h\ot k)\varepsilon(l)
\end{equation}
for all $h,k,l\in H$. \eqref{3.3.1.1} means in fact that $\theta$ is right
$H$-linear. In terms of the corresponding integral maps $\gamma\in V_4$,
\eqref{3.3.1.1} can be rewritten as follows:
\begin{equation}\eqlabel{3.3.1.2}
h\vdash \gamma=\gamma \dashv h
\end{equation}
for all $h\in H$, where the actions of $H$ on $\Hom(H,H^*)$ are given by
\begin{eqnarray*}
(h\vdash\gamma)(g)(k)&=& \sum h_{(3)}\gamma(gh_{(1)})(kh_{(2)})\\
(\gamma\dashv h)(g)(k)&=& \gamma(g)(k)h
\end{eqnarray*}
for all $h,g,k\in H$. Our results can be summarized as follows.

\begin{theorem}\thlabel{3.3.1}
Let $H$ be a flat Hopf algebra over $k$, and consider the
Doi-Hopf datum $(k,k,H)$. There exist $k$-linear maps
$$i:\ \int_{H^*}^l\to V_3\cong V_4~~{\rm and}~~
p:\ V_4\cong V_3\to \int_{H^*}^l$$
such that $p\circ i=I_{\int_{H^*}^l}$. Furthermore, the image of $i$ in
$V_3$ (resp. $V_4$) is the submodule consisting of maps that satisfy
\eqref{3.3.1.1} (resp. \eqref{3.3.1.2}).
\end{theorem}

\subsection{The finite case}\selabel{3.4}
Let $(H,C,A)$ be a Doi-Hopf datum such that
$H$ is a Hopf algebra and
$C$ is finitely generated and projective over $k$. Then
${}^C{\cal M}(H)_A\cong {\cal M}_{A\#C^*}$ and the forgetful functor
$F:\ {}^C{\cal M}(H)_A\to {\cal M}_A$ is isomorphic to the restriction of
scalars
functor $R:\ {\cal M}_{A\#C^*}\to {\cal M}_A$.
Moreover, $i:\ A\#C^*\to \Hom(C,A)$ is an
isomorphism of algebras and of $(A,A\#C^*)$-bimodules.
Therefore an $A$-integral $\gamma$ can be viewed
as a right $C^*$-module map $\gamma:\ C\to A\#C^*$ which is
centralized by the action of $A$ on $\Hom(C, A\#C^*)$.
Recall that $A\#C^*$ is a right
$C^*$-module after restriction of scalars via the map
 $i_{C^*}:\ C^*\to A\#C^*$.
Identifying $A\#C^*\cong\Hom(C,A)$, we find that
$\Hom(C, A\#C^*)$ is an $A$-bimodule (see (\ref{eq:1.11a},\ref{eq:1.11b})).
The left and right action of $A$ on $\Hom(C, A\#C^*)$ are given by
\begin{eqnarray*}
(a\rhu \gamma)(c)&=&\sum (a_{<0>}\# \varepsilon) \gamma (c \cdot a_{<-1>})\\
(\gamma \lhu a)(c)&=&\gamma (c) (a\#\varepsilon)
\end{eqnarray*}
for all $a\in A$, $\gamma \in \Hom(C,A\#C^*))$ and $c\in C$.\\
>From \thref{2.1.3} we obtain the following necessary and sufficient condition
for the separability of the extension $A\to A\#C^*$.

\begin{corollary}\colabel{3.4.1}
Let $(H,C,A)$ be a Doi-Hopf datum, with $H$ a Hopf algebra, and $C$ finitely
generated and projective as a $k$-module.
The following statements are equivalent:\\
1) the extension $A\to A\#C^*$ is separable;\\
2) there exists  a normalized $A$-integral $\gamma:\ C\to A\#C^*$.\\
In this situation,
if $A$ is semisimple artinian, then $A\# C^*$ is semisimple artinian also.
\end{corollary}

\begin{remarks}\relabel{3.4.2}\rm
1) In Corollary 3.3 of \cite{CMZ3} we proved that if there exists
$\gamma :C\to A\#C^*$ a normalized $A$-integral then the extension $A\to
A\#C^*$ is right
semisimple. As any separable extension is a semisimple extension,
the above Corollary improves
Corollary 3.3 of \cite{CMZ3} and also gives us the converse.\\
2) Consider the particular case $C=H$. Then $A$ is an
$H^*$-module algebra and we get a necessary and sufficient condition for the
extension  $A\to A\#H^*$ to be separable.\\
Several results connected to the separability of this extension
$A\to A\#H^*$ have appeared in the literature.\\
In the first place, if $H$ is finitely dimensional
and cosemisimple, then the extension $A\to A\#H^*$ is separable
(see \cite[Theorem 4]{CF1}).\\
Secondly, Propositions 1.3 and 1.5 of \cite{VXZ} give necessary and
sufficient condition for the  extension $A\to A\#H^*$ to be separable in
the cases
where $A$ is an $H$-Galois extension and
$H$ contains a cocommutative integral.\\
Finally, Theorem 3.14 of \cite{DT} can also be connected to our results:
$A\#H^*$ is a right $H^*$-comodule algebra via $I_A\ot\delta_{H^*}$, and
this makes the extension $A\to A\#H^*$ an $H^*$-Galois extension. Now
$$(A\#H^*)^A:=\{a\#h^* \mid (b\#\varepsilon)(a\#h^*)=(a\#h^*)(b\#\varepsilon)~
{\rm for~all~}b\in A \}$$
is a right $H^*$-module algebra, by the Miyashita-Ulbrich action, and
therefore a left $H$-comodule algebra, since $H$ is finitely generated
projective. From \cite[Theorem 3.14]{DT},
we obtain that the extension $A\to A\# H^*$ is separable if and only if
there exists  a total integral $\varphi:\ H\to (A\#H^*)^A$.
Disadvantages of this approach are the lack of control that we have
on the space $(A\#H^*)^A$, and also the fact that the $H$-coaction
coming from the Miyashita-Ulbrich action is not very handable.\\
3) Let $C=H$, and assume that $H$ is finite dimensional over a field $k$.
Let $t$
be a right integral in $H$ and $\Lambda$ a right integral in $H^*$. Then $H$
is free and cyclic as a right $H^*$-module with basis $\{t\}$:
$H=t\lhua H^*$ (see \cite{Rd1} or \cite{S}). The inverse
of the map
$$H^*\to H:\ h^*\mapsto t\lhua h^*$$  is the map
$$H\to H^*:\ h\mapsto S^{-1}(h)\rhua \Lambda$$
Thus right $H^*$-linear maps
$\gamma:\ H\to A\#H^*$ correspond to elements
$s=\sum_i a_i\ot h_{i}^{*}\in A\ot H^*$. From
\coref{3.4.1}, it follows that the separability of the extension $A\to
A\#H^*$ is
equivalent to the existence of
an element $s=\sum_i a_i\ot h_{i}^{*}\in A\ot H^*$
satisfing the two following conditions.
$$\sum a_{<0>} a_{i_{<0>}}\ot \lan h_{i}^{*}, l_{(1)}a_{<-2>}\ran
l_{(2)}a_{<-1>}a_{i_{<-1>}}S^{-1}(a_{<-3>}) =
\sum a_{i_{<0>}}a \ot \lan h_{i}^{*}, l_{(1)}\ran l_{(2)}a_{i_{<-1>}}$$
$$\sum \lan h_{i}^{*}, h_{(2)}\ran
\lan S^{-1}(h_{(1)}) \rhua \Lambda , h_{(3)}a_{i_{<-1>}}\ran a_{i_{<0>}}=
\varepsilon (h) 1_A$$
for all $a\in A$, $l$, $h\in H$.
\end{remarks}

Now take $H=A=C$ and consider the Doi-Hopf datum $(H,H,H)$, with
$H$ a finitely generated and projective Hopf algebra. The smash product
${\cal H}(H)=H\#H^*$ is usually called the {\sl Heisenberg double} of $H$.
Applying \coref{3.4.1}, we find the following necessary and sufficient
conditions for $H\to {\cal H}(H)$ to be separable.
\begin{corollary}\colabel{3.4.3}
Let $H$ be a finitely generated and projective Hopf algebra over $k$. The
following statements are equivalent:\\
1) the extension $H\to {\cal H}(H)$ is separable;\\
2) there exists  a normalized $H$-integral $\gamma:\ H\to {\cal H}(H)$.
\end{corollary}

\subsection{Yetter-Drinfel'd modules}\selabel{3.5}
In this section $H$ will be a $k$-flat Hopf algebra with bijective antipode
$S$.
Recall that a right-left Yetter Drinfel'd module $M$ is a $k$-module, which is
at once a left $H$-comodule and a right $H$-module, such that the following
compatibility relation holds:
\begin{equation}\eqlabel{3.5.1.1}
\sum  m_{<-1>}h_{(1)}\otimes m_{<0>}h_{(2)} =
\sum S^{-1}(h_{(3)})m_{<-1>}h_{<1>}\otimes m_{<0>}h_{(2)}
\end{equation}
for all $h\in H$ and $m\in M$.
The category of right-left Yetter-Drinfel'd
modules and $H$-linear $H$-colinear maps will be denoted by
${}^H{\cal YD}_H$. \\
In \cite{CMZ1} it is shown that there is a category isomorphism
$${}^H{\cal YD}_H\cong {}^H{\cal M}(H\otimes H^{op})_H$$
The left $H\otimes H^{\rm op}$-coaction and right $H\otimes H^{\rm op}$-action
on $H$ are given by the formulas
\begin{eqnarray}
\rho(h)&=&\sum h_{(1)}\otimes S^{-1}(h_{(3)})\otimes h_{(2)}\eqlabel{3.5.1.2}\\
l\cdot (h\otimes k)&=&klh \eqlabel{3.5.1.3}
\end{eqnarray}
for all $h$, $k$, $l\in H$.
>From now on, we will identify the categories ${}^H{\cal YD}_H$ and
${}^H{\cal M}(H\otimes H^{op})_H$ using the isomorphism from \cite{CMZ1}.\\
We recall some basic facts about the Drinfel'd double, as introduced in
\cite{Dr}. Our main references are \cite{Ma} and \cite{Rd}.
Let $H$ be a finitely generated projective Hopf algebra.
Then the antipode $S_H$ is bijective (see \cite{P}), and
the Drinfel'd double $D(H)=H\bowtie H^{*\,\rm cop}$ is defined as follows:
$D(H)=H\otimes H^{*}$ as a $k$-module, with multiplication,
comultiplication, counit and antipode given
by the formulas

\begin{eqnarray*}
(h\bowtie f)(h^{\prime}\bowtie f^{\prime})&=&\sum h_{(2)}h^{\prime}
\bowtie f*\lan f^{\prime},S^{-1}h_{(3)}?h_{(1)}\ran\\
\Delta_{D(H)}(h\bowtie f)&=&\sum (h_{(1)}\bowtie f_{(2)})\otimes
(h_{(2)}\bowtie f_{(1)})\\
\varepsilon_{D(H)}&=&\varepsilon_H\otimes\varepsilon_{H^{*\,\rm cop}}\\
S_{D(H)}(h\bowtie f)&=&\sum f_{(2)}\rhua (S_Hh_{(1)})\bowtie
(S(h_{(2)})\rhua (S_{H^{*}}f_{(1)})
\end{eqnarray*}
for $h,h^{\prime}\in H$ and $f,f^{\prime}\in H^{*}$.\\
Consider the Doi-Hopf datum $(H\ot H^{op}, H, H)$, with structures given by
(\ref{eq:3.5.1.2}-\ref{eq:3.5.1.3}). $H^*$ is a left
$H\otimes H^{\rm op}$-module algebra, with
$$\lan (h\otimes k)\trr h^{*},l\ran =\lan h^{*},klh\ran $$
for all $h,k,l\in H$ and $h^{*}\in H^{*}$,
and in \cite{CMZ1}, it is shown that the Drinfel'd double is the smash
product
$$D(H)=H\#H^{*}$$
For a Hopf algebra $H$ that is not necessary finitely generated projective,
consider Koppinen's smash product ${\cal D}(H)=\#(H,H)$. If $H$ is
finitely generated and projective, then ${\cal D}(H)\cong D(H)$
(see \cite{Tak} for a similar result, where it was proved that
$D(H)\cong \End(H)^*$).
We can therefore view ${\cal D}(H)$ as a generalization of the Drinfel'd
double to
the case of infinite dimensional Hopf algebras. The structure of
${\cal D}(H)$ is the following: as a $k$-module, ${\cal D}(H)=\End(H)$,
and the multiplication is given by the formula
$$(f\bu g)(h)=\sum
f(h_{(1)})_{(2)} g\Bigl (S^{-1}(f(h_{(1)})_{(3)}) h_{(2)}
f(h_{(1)})_{(1)}\Bigl)$$
${\cal D}(H)$ is a right $H^*$-module via
$$(f\cdot h^*)(h)=\sum
f(h_{(1)})_{(2)} \lan h^*,  S^{-1}\Bigl(f(h_{(1)})_{(3)}\Bigl) h_{(2)}
f(h_{(1)})_{(1)}\ran$$
and $\Hom (H, {\cal D}(H))$ is an $H$-bimodule via the formulas
$$(g\rhu \gamma)(h)(l)=\sum g_{(3)}\gamma \Bigl (S^{-1}(g_{(5)})hg_{(1)} \Bigr)
\Bigl(S^{-1}(g_{(4)})lg_{(2)}\Bigr)$$
and
$$(\gamma \lhu g)(h)(l)=\gamma (h)(l)g$$
for any $g$, $h$, $l\in H$, and $\gamma :H\to {\cal D}(H)$.

We can now apply the Maschke \thref{2.1.5} to the Doi-Hopf datum
$(H\ot H^{op}, H, H)$. An $A=H$-integral for $(H\ot H^{op}, H, H)$
will be called a {\it quantum $H$-integral}. This can be reformulated
as follows.

\begin{definition}\delabel{3..5.1}
Let $H$ be a Hopf algebra with a bijective antipode. A quantum $H$-integral
is a right $H^*$-module map $\gamma:\ H\to {\cal D}(H)$ which is centralized
by the left and right $H$-action, that is
$$g\rhu \gamma=\gamma \lhu g$$
for all $g\in H$. $\gamma$ is a total quantum $H$-integral if
$$\sum \gamma(h_{(1)})h_{(2)}=\varepsilon (h) 1_H$$
for all $h\in H$.
\end{definition}

>From \thref{2.1.5} we obtain immediately the following version of Maschke's
Theorem for Yetter-Drinfel'd modules:

\begin{corollary}\colabel{3.5.2}
Let $H$ be a Hopf algebra. Assume that the antipode is bijective,
and that $H$ is projective as a $k$-module. Then
the following statements are equivalent:\\
1) the forgetful functor $F:\ {}^H{\cal YD}_H\to {\cal M}_H$ is separable;\\
2) there exists a total quantum $H$-integral $\gamma:\ H\to {\cal D}(H)$.
\end{corollary}

If $H$ is finite, then we obtain the following result:

\begin{corollary}\colabel{3.5.3}
Let $H$ be a finitely generated and projective Hopf algebra.
The following statements are equivcalent:\\
1) The extension $H\to D(H)$ is separable;\\
2) there exists a total quantum $H$-integral $\gamma:\ H\to {\cal D}(H)
\cong D(H)$.
\end{corollary}

Let $H$ be a finite dimensional Hopf algebra over a field of characteristic
zero. A classical result of Larson and Radford \cite{LR} tells us that
$H$ is semisimple if and only if $H$ is cosemisimple. We will now show
that the space of quantum $H$-integrals measures how far semisimple
Hopf algebras are from cosemisimple Hopf algebras, if we work over an
arbitrary field. First recall from \cite{Rd} that $D(H)$ is semisimple if and
only if $H$ is semisimple and cosemisimple.

\begin{corollary}\colabel{3.5.4}
Let $H$ be a finite dimensional semisimple Hopf algebra over a field $k$.
The following statements are equivalent:\\
1) $H$ is cosemisimple;\\
2) there exists a total quantum $H$-integral $\gamma:\ H\to {\cal D}(H)
\cong D(H)$.
\end{corollary}

\begin{proof}
If $H$ is cosemisimple, then $D(H)$ is semisimple, and therefore a
separable extension of the field $k$. But then
$D(H)$ is separable over $H$ (see for example \cite[Lemma 1.1]{NvBvO},
and the first implication follows from \coref{3.5.3}.\\
Conversely, if there exists total quantum $H$-integral, then
$H\to D(H)$ is separable, and it follows that $D(H)$ is semisimple, since
$H$ is semisimple. From Radford's results (\cite{Rd}), it follows that
$H$ is cosemisimple.
\end{proof}

\begin{remark}\rm
If $H$ is finite dimensional, then $H$ is free of rank one as a right
$H^*$-module, generated by a nonzero right integral $t\in H$.
The existence of a total quantum $A$-integral $\gamma:\ H\to D(H)$ is then
equivalent to the existence of an element $z(= \gamma(t))\in D(H)$
satisfying two properties. These properties are the translations of the
facts that $\gamma$ is centralized by the left and right action of $H$
and that $\gamma$ is normalized. Unfortunately, these properties cannot
be written down in an elegant and transparent way..
\end{remark}

\subsection{Long's category of dimodules}\selabel{3.7}
Consider the Doi-Hopf datum $(H,A=H,C=H)$, where,
$A=H$ is a left $H$-comodule algebra via $\Delta$ and $C=H$ is a
right $H$-module coalgebra with the trivial $H$-action, that is
$g\cdot h=g\varepsilon (h)$, for all $g$, $h\in H$. The compatibility
relation for Doi-Hopf modules now takes the form
$$\rho(mh)=\sum m_{<-1>} \ot m_{<0>}h$$
for all $m\in M$ and $h\in H$. If $H$ is commutative and cocommutative,
then we obtain $H$-{\it dimodules} in the sense of Long \cite{L}, and
in this situation, dimodules coincide with Yetter-Drinfel'd modules.
Dimodules and dimodule algebras have been a basic tool in the study of
the generalizations of the Brauer-Wall group, and, recently,
the third author investigated the relation between the category of
dimodules and certain nonlinear equation (see \cite{Mi2}).\\
The category of Long dimodules will be denoted by ${}^H{\cal L}_H$.
If $H$ is finitely generated and projective, then
${}^H{\cal L}_H\cong {\cal M}_{H\ot H^*}$. For a group $G$, a $kG$-dimodule
is a $G$-graded module $M$, on which the group $G$ acts in such a way that
the homogeneous components are themselves $G$-modules:
$$g\cdot M_h\subset M_h$$
for all $g,h\in G$.\\
It is easy to see that the multiplication on the Koppinen smash product
$\#(H,H)=\End(H)$ is nothing else then the convolution
$$(f\bu g)(h)=(f * g)(h)=\sum f(h_{(1)}) g(h_{(2)})$$
for any $f$, $g\in \End(H)$ and $h\in H$.
Furthermore, the right $H^*$-action on $\#(H,H)=\End(H)$ is given by
$$(f\cdot h^*)(h)=\sum f(h_{(1)}) \lan h^*, h_{(2)}\ran =f(h^*\rhua h)$$
for all $f\in \End(H)$, $h^*\in H^*$, and $h\in H$, and the $H$-bimodule
structure of $\Hom (H, \End(H))$ is trivial:
$$(g\rhu \gamma)(h)(l)=g\gamma (h)(l)~~{\rm and}~~
(\gamma \lhu g)(h)(l)=\gamma (h)(l)g$$
Applying \thref{2.1.5}, we now find a Maschke Theorem for dimodules:
\begin{corollary}\colabel{3.7.1}
For a projective Hopf algebra over a commutative ring $k$, the
following statements are equivcalent:\\
1) the forgetful functor $F:\ {}^H{\cal L}_H\to {\cal M}_H$ is separable;\\
2) there exists an $H^*$-linear map $\gamma:\ H\to \End(H)$
such that $\im(\gamma(h))\subset Z(H)$ and
$\sum \gamma (h_{(1)})(h_{(2)})$ = $\varepsilon (h)1_H$, for all $h\in H$.
\end{corollary}

\subsection{Modules graded by $G$-sets}\selabel{3.8}
Let $G$ be a group, $X$ a right $G$-set and $H=kG$ a group algebra.
Then the grouplike coalgebra $C=kX$
is a right $kG$-module coalgebra, and a left $kG$-comodule
algebra $A$ is nothing else than a $G$-graded $k$-algebra (see
for example \cite{M}). In this case the category of Doi-Hopf
module ${}^{kX}{\cal M}(kG)_{A}$ is the category gr-$(G,X,A)$ of right
$X$-graded $A$-modules (see \cite{NRvO}). An object in this category is a right
$A$-module such that  $M=\oplus_{x\in X} M_x$ and
$M_xA_{g}\subseteq M_{x\cdot g}$, for all
$g\in G$ and $x\in X$.\\
Applying \thref{2.1.5}, we obtain the following result. Observe that,
in the case where $X=G$, the
second statement of the next Theorem can be found in \cite{NvBvO}.

\begin{proposition}\prlabel{3.8.1}
Let $G$ be a group, $X$ a right $G$-set and $A$ be a $G$-graded $k$-algebra.
Then the map $\gamma:\ kX\to \Hom(kX, A)$ given by the formula
$$\gamma(x)(y)=\delta_{x,y} 1_A$$
for all $x$, $y\in X$, is a normalized $A$-integral.\\
Consequently, the forgetful functor
$F:\ {\rm gr}\hbox{-}(G,X,A)\to {\cal M}_A$ is
separable.
\end{proposition}

\begin{proof}
Let $\{ x, p_x~|~x\in X \}$ be the canonical dual basis of $kX$. For
all $x$, $y$, $z\in X$ we have
$$\gamma (x\lhua p_y)(z)=\gamma (\lan p_y, x\ran x)(z)=
\delta_{x,y}\delta_{x,z}1_A$$
and
\begin{eqnarray*}
(\gamma(x)\cdot p_y)(z)&=&\sum \gamma (x)(z)_{<0>}
\lan p_y, z\cdot \gamma (x)(z)_{<-1>} \ran\\
&=&\delta_{x,z}1_A \lan p_y, z \ran\\
&=&\delta_{x,z}\delta_{y,z}1_A
\end{eqnarray*}
and this shows that $\gamma$ is right $(kX)^*$-linear. Now, let $g\in G$ and
$a_g\in A_g$ (this means that $\rho_A(a_g)=g\ot a_g$).
Then
\begin{eqnarray*}
(a_g \rhu \gamma)(x)(y)&=& (a_g\cdot \gamma (x\cdot g))(y)\\
&=& a_g\gamma (x\cdot g)(y\cdot g)\\
&=& \delta_{x\cdot g, y\cdot g}a_g\\
&=& \delta_{x,y}a_g
\end{eqnarray*}
and
$$(\gamma \lhu a_g)(x)(y)=(\gamma (x)\cdot a_g)(y)=
\gamma (x)(y)a_g=\delta_{x,y}a_g$$
and $\gamma$ is centralized by the action of $A$. Hence
$\gamma$ is a normalized $A$-integral.
\end{proof}

We now focus attention to the right adjoint
$$G=kX\ot\bullet:\ {\cal M}_A\to {\rm gr}\hbox{-}(G,X,A)$$
of the forgetful functor $F$. This case is interesting because,
in the particular situation where $X=G$ and $A$ is a strongly
$G$-graded algebra, the separability of the functor $G$ is equivalent
to the separability of $A$ as an $A_1$-algebra.\\
>From \thref{2.2.1}, we obtain the following result, which was already
shown using other methods in \cite[Theorem 3.6]{NvBvO} (if $X=G$) and
in \cite[Sec. 4]{R} (for general $X$).

\begin{corollary}\colabel{3.8.2}
Let $G$ be a group, $X$ a right $G$-set and $A$ a $G$-graded $k$-algebra.
The following statements are equivalent:\\
1) the functor $G=kX\ot \bu:\ {\cal M}_A\to {\rm gr}\hbox{-}(G,X,A)$
is separable;\\
2) there exists a normalized dual $A$-integral $z=\sum_i x_i\ot a_i \in
kX\ot A$,
that is
$$\sum_i a_i=1_A$$
and
$$\sum x_i\ot b_g a_i=\sum x_i\cdot g\ot a_ib_g$$
for all $g\in G$ and $b_g\in A_g$.
\end{corollary}

We will further investigate the second condition of \coref{3.8.2}. First
recall that a $G$-subset $X'$ of $X$ is a subset $X'\subset X$ such that
 $x'\cdot g\in X'$, for all $x'\in X'$, $g\in G$.

\begin{corollary}\colabel{3.8.3}
Let $G$ be a group, $X$ a right $G$-set and $A$ be a $G$-graded $K$-algebra.\\
1) Let $X'$ be a finite $G$-subset of $X$
such that $\#(X')$ is invertible in $k$.
Then the functor $G=kX\ot \bu:\ {\cal M}_A\to {\rm gr}\hbox{-}(G,X,A)$
is separable.\\
2) If $A$ is strongly graded, and the functor
$G:\ {\cal M}_A\to {\rm gr}\hbox{-}(G,X,A)$ is separable, then there exists
a finite
$G$-subset $X'$ of $X$.
\end{corollary}

\begin{proof}
1) We have a dual total $A$-integral of $kX\ot A$, namely
$$z={1\over \#(X')} \sum_{x'\in X'} x'\ot 1_A$$
2) Let $0\neq z=\sum_{i=1}^n x_i\ot a_i$, with $x_i\in X$ and $a_i\in A$, be
 dual integral. We take $n$ as small as possible. Then $a_i\neq 0$, for any
$1\leq i\leq n$. We claim that $\{x_1,\ldots,x_n\}$ is a $G$-subset of $X$.
It suffices to show that $x_i\cdot g\in \{x_1,\ldots,x_n\}$, for any
$i\in\{1,\ldots, n\}$ and $g\in G$.\\
Since $A$ is strongly graded, we can find
finite sets $\{b_1,\ldots,b_m\}\subset A_g,~
\{b'_1,\ldots,b'_m\}\subset A_{g^{-1}}$ such that
$$\sum_{j=1}^m b_jb'_j=1$$
For all $j$, we have, using the fact that $z$ is a dual integral:
$$\sum_{i=1}^n x_i\ot b_ja_i=\sum_{i=1}^n x_i\cdot g\ot a_ib_j$$
and
$$\sum_{i=1}^n x_i\cdot g\ot a_i=
\sum_{i=1}^n\sum_{j=1}^m x_i\cdot g\ot (a_ib_jb'_j)=
\sum_{i=1}^n x_i\ot \sum_{j=1}^m b_ja_ib'_j$$
and the statement follows from the fact that $X$ is a basis of $kX\ot A$
as a right $A$-module, and the fact that $a_i\neq 0$.
\end{proof}

In particular, we have the following result for a field of characteristic
zero.

\begin{corollary}\colabel{3.8.4}
Let $k$ be a field of characteristic zero, $G$ a group, $X$ a right $G$-set,
and $A$ a strongly $G$-graded $k$-algebra.
The following statements are equivalent:\\
1) The functor $G=kX\ot \bu:\ {\cal M}_A\to {\rm gr}\hbox{-}(G,X,A)$
is separable.\\
2) there exists a finite $G$-subset $X'$ of $X$;\\
3) There exists $x\in X$ with finite $G$-orbit ${\cal O}(x)$.
\end{corollary}

\begin{remarks}\relabel{3.8.5}\rm
1) In the second part of \coref{3.8.3}, we need the assumption that $A$ is
strongly graded. An easy counterexample is the following: let $X$ be any
$G$-set without finite orbit, and $A$ an arbitrary $k$-algebra. We give
$A$ the trivial $G$-grading. Then for any $x\in X$, $x\ot 1_A$ is a normalized
dual integral.\\
2) Let $X=G=\nint$, and $A$ a strongly $G$-graded $k$-algebra. Then the functor
$G:\ {\cal M}_A\to {\rm gr}\hbox{-}A$ is not separable.\\
3) Let $X=G$, where $G$ acts on $G$ by conjugation, that is
$g\cdot h=h^{-1}gh$, for all $g$, $h\in G$. The corresponding Doi-Hopf modules
will be called {\it right $G$-crossed $A$-modules}, and we will call
denote the category of right $G$-crossed $A$-modules by
${\rm gr}\hbox{-}(G,G,A)={}^G{\cal C}_A$.
$z=1_G\ot 1_A$ is
a dual normalized $A$-integral, hence the functor
$F:\ {\cal M}_A\to {}^G{\cal C}_A$ is separable. If $A=kG$, then the category
${}^G{\cal C}_A$ is just the category of crossed $G$-modules defined
by Whitehead (see \cite{W}).
\end{remarks}

Now take $X=G$, where the action is the usual multiplication of $G$. Then
the category gr-$(G,G,A)$, also denoted by gr-$A$ is the category of
$G$-graded $A$-modules, and we obtain the following properties:\\
1) if $F:\ {\cal M}_A\to {\rm gr}\hbox{-}A$ is separable then $G$ is finite;\\
2) if $G$ is finite and $\# G$ is invertible in $k$, then $F$
is a separable functor.\\
In particular, for a strongly graded $G$-algebra $A$ (i.e.. a $kG$-Hopf
Galois extension of $A_1$), we obtain a necessary and sufficient condition
for the extension $A_1\subset A$ to be separable, using \thref{3.1.5}
and \coref{3.8.2}. Compare this to \cite[Proposition 2.1]{NvBvO}.\\
If $A$ is a strongly graded $G$-algebra with $G$ a finite group and
$\# G$ invertible in $k$, then $A$ is a separable extension of $A_1$.\\
In the next example we will construct a strongly graded $G$-algebra $A$ such
that
$A/A_1$ is a separable extension, with  $\# G$ not invertible in $k$.

\begin{example}\exlabel{3.8.6}\rm
Let $k=\nint_2$ and $C_2=\{1,c\}$ the cyclic group of order two. On
the matrix ring $A={\rm M}_2(\nint_2)$, we consider the following grading
(see \cite{DINR}): $A_1$ and $A_c$ are the vector spaces with respectively
$$\left\{\pmatrix{1&0\cr 0&1\cr},\pmatrix{0&1\cr 1&1\cr}\right\}~~{\rm and}~~
\left\{\pmatrix{0&1\cr 1&0\cr},\pmatrix{1&1\cr 0&1\cr}\right\}$$
as $k$-basis.
Then $A$ is a strongly $G$-graded ring (even a crossed
product), and
$$z=1\ot \pmatrix{0&1\cr 1&1\cr}+c\ot \pmatrix{1&1\cr 1&0\cr}$$
is a normalized dual $A$-integral (this can be proved by an easy computation on
the basis elements). Hence $A$ is a separable extension of $A_1$, and,
obviously, $\# C_2=2={\rm char} (k)$.
\end{example}

We end this section pointing out that, over a field $k$ of characteristic zero,
and for a strongly $G$-graded $k$-algebra $A$, the extension $A_1\subset A$
is separable if and only if $G$ is a finite group.

{\bf Acknowledgement} The authors thank Paul Taylor for his kind permission to
use the "diagrams" software.

\end{document}